\documentclass{article}
\usepackage{amsmath,amssymb,amsthm,graphicx,epsfig,subfigure,float,url}
\usepackage[colorlinks=true]{hyperref}
\usepackage{pdfsync}

\usepackage{verbatim}
\usepackage[applemac]{inputenc}

\def\R{\mathbb{R}}

\def\Z{\mathbb{Z}}
\def\C{\mathbb{C}}

\renewcommand{\geq}{\geqslant}
\renewcommand{\leq}{\leqslant}

\newtheorem{theorem}{Theorem}

\newtheorem{definition}{Definition}
\newtheorem{lemma}{Lemma}
\theoremstyle{definition}
\theoremstyle{definition}\newtheorem{remark}{Remark}

\usepackage{graphics} 
\usepackage{epsfig} 
\usepackage{epstopdf}
\usepackage{subfigure}
\usepackage{color}

\graphicspath{{./Figures/}}
\allowdisplaybreaks[4]

\title{\bf PI regulation control of a 1-D semilinear wave equation}

\author{Hugo Lhachemi\footnote{University College Dublin, Dublin, Ireland (\texttt{hugo.lhachemi@ucd.ie}).}
\and
Christophe Prieur\footnote{Universit\'e Grenoble Alpes, CNRS, Grenoble-INP, GIPSA-lab, F-38000, Grenoble (\texttt{christophe.prieur@gipsalab.fr}).}
\and
Emmanuel Tr\'elat\footnote{Sorbonne Universit\'e, CNRS, Universit?e de Paris, Inria, Laboratoire Jacques-Louis Lions (LJLL), F-75005 Paris, France (\texttt{emmanuel.trelat@sorbonne-universite.fr}).}
}

\begin{document}
\maketitle

This paper is concerned with the Proportional Integral (PI) regulation control of the left Neumann trace of a one-dimensional semilinear wave equation. The control input is selected as the right Neumann trace. The control design goes as follows. First, a preliminary (classical) velocity feedback is applied in order to shift all but a finite number of the eivenvalues of the underlying unbounded operator into the open left half-plane. We then leverage on the projection of the system trajectories into an adequate Riesz basis to obtain a truncated model of the system capturing the remaining unstable modes. Local stability of the resulting closed-loop infinite-dimensional system composed of the semilinear wave equation, the preliminary velocity feedback, and the PI controller, is obtained through the study of an adequate Lyapunov function. Finally, an estimate assessing the set point tracking performance of the left Neumann trace is derived.

\section{Introduction}

Due to its widespread adoption by industry~\cite{aastrom1995pid,astrom2008feedback}, the stabilization and regulation control of finite-dimensional systems by means of Proportional-Integral (PI) controllers has been intensively studied. For this reason, the opportunity of extending PI control strategies to infinite-dimensional systems, and in particular to systems modeled by partial differential equations (PDEs), has attracted much attention in the recent years. Efforts in this research direction were originally devoted to the case of bounded control operators~\cite{pohjolainen1982robust,pohjolainen1985robust} and then extended to unbounded control operators~\cite{xu1995robust}. The study of PI control design combined with high-gain conditions was reported in~\cite{logemann1992robust}. More recently, the problem of PI boundary control of linear hyperbolic systems has been reported in a number of works~\cite{bastin2015stability,dos2008boundary,lamare2015control,xu2014multivariable}. This research direction has then been extended to the case of nonlinear transport equations~\cite{bastin2019exponential,coron2019pi,martins2014design,hayat2019pi,rodrigues2013lmi,trinh2017design}. The case of the boundary regulation control of the Neumann trace for a linear reaction-diffusion in the presence of an input delay was considered in~\cite{lhachemi2019pi}. The case of the boundary regulation control of the boundary velocity for linear damped wave equations, in the presence of a nonlinearity in the boundary conditions, has been considered in~\cite{barreau2019practical,terrand2018regulation}. A general procedure allowing the addition of an integral component for regulation control to open-loop exponentially stable semigroups with unbounded control operators has been proposed in~\cite{terrand2018lyapunov,terrand2019adding}. 

This paper is concerned with the PI regulation control of the left Neumann trace of a one-dimensional semilinear (undamped) wave equation. The selected control input takes the form of the right Neumann trace. The control design procedure goes as follows. First, inspired by~\cite{coron2006global}, a preliminary (classical) velocity-feedback is applied in order to shift all but a finite number of the eigenvalues of the underlying unbounded operator into the open left half-plan. Then, inspired by the early work~\cite{russell1978controllability} later extended in~\cite{coron2004global,coron2006global,schmidt2006} to semilinear heat and wave PDEs, we leverage on the projection of the system trajectories into a Riesz basis formed by the generalized eigenstructures of the unbounded operator in order to obtain a truncated model capturing the remaining unstable modes. Finally, similarly to~\cite{lhachemi2019pi}, this finite dimensional model is augmented to include the integral component of the PI controller, allowing to compute a stabilizing feedback. The local stability of the resulting closed-loop infinite-dimensional system, and the subsequent set point regulation performance, is assessed by a Lyapunov-based argument. The theoretical results are illustrated based on the simulation of an open-loop unstable semilinear wave equation.

The paper is organized as follows. The investigated control problem is introduced in Section~\ref{sec: problem setting}. The proposed control design procedure is presented in Section~\ref{sec: control design} in a comprehensive manner. The subsequent stability analysis is carried out in Section~\ref{sec: stability and regulation} while the theoretical results are numerically illustrated in Section~\ref{sec: numerical illustration}. Finally, concluding remarks are formulated in Section~\ref{sec: conclusion}.

\section{Problem setting}\label{sec: problem setting}

Let $L>0$ and let $f : \R \rightarrow \R$ be a function of class $\mathcal{C}^2$. We consider the following wave equation on $(0,L)$:
\begin{subequations}\label{eq: wave equation}
\begin{align}
&\dfrac{\partial^2 y}{\partial t^2} = \dfrac{\partial^2 y}{\partial x^2} + f(y) , \\
& y(t,0) = 0 , \qquad \dfrac{\partial y}{ \partial x}(t,L) = u(t) ,\\
& y(0,x) = y_0(x) , \qquad \dfrac{\partial y}{\partial t}(0,x) = y_1(x) , 
\end{align}
\end{subequations}
for $t > 0$ and $x \in (0,L)$, where the state is $y(t,\cdot):[0,L]\rightarrow\R$ and the control input is $u(t)$ and applies to the right Neumann trace. The control objective is to design a PI controller in order to locally stabilize the closed-loop system and locally regulate the system output selected as the left Neumann trace:
\begin{equation}\label{eq: system output}
z(t) = \dfrac{\partial y}{ \partial x}(t,0) .
\end{equation}

\begin{definition}\label{def: steady state}
A function $y_e \in\mathcal{C}^2([0,L])$ is a steady-state of (\ref{eq: wave equation}) with associated constant control input $u_e \in \R$ and constant system output $z_e \in \R$ if
\begin{equation*}
\begin{split}
& \dfrac{\mathrm{d}^2 y_e}{\mathrm{d} x^2}(x) + f(y_e(x)) = 0 , \qquad x \in (0,L), \\
& y_e(0) = 0 , \quad \dfrac{\mathrm{d} y_e}{\mathrm{d} x}(L) = u_e , \\
& z_e = \dfrac{\mathrm{d} y_e}{\mathrm{d} x}(0) .
\end{split}
\end{equation*}
\end{definition}

\begin{remark}\label{rem: existence of steady states}
Introducing $F(y) = \int_0^y f(s) \,\mathrm{d}s$ for any $y \in \R$, assume that one of the two following properties holds:
\begin{itemize}
\item $F(y) \rightarrow + \infty$ when $\vert y \vert \rightarrow +\infty$;
\item for any $a > 0$, when it makes sense, the integral $\int \frac{\mathrm{d}y}{\sqrt{a-F(y)}}$ diverges at $-\infty$ and $+\infty$.
\end{itemize}
Then we have the existence of a steady state $y_e \in\mathcal{C}^2([0,L])$ of (\ref{eq: wave equation}) associated with any given value of the system output $z_e \in \R$. Indeed, define $y \in \mathcal{C}^2([0,l))$ with $0 < l \leq + \infty$ as the maximal solution of $y''+f(y)=0$ with $y(0)=0$ and $y'(0) = z_e$. We only need to assess that $l > L$. Multiplying by $y'$ both sides of the ODE satisfied by $y$ and then integrating over $[0,x]$, we observe that $y$ satisfies the conservation law $y'(x)^2 + 2F(y(x)) = z_e^2$ for all $x \in [0,l)$. Hence any of the two above assumptions implies that $y$ and $y'$ are bounded on $[0,l)$. Thus $l = + \infty$ and the associated steady state control input is given by $u_e = y'(L)$.
\end{remark}

Given a desired value of the system output $z_e \in\R$ and an associated steady state function $y_e \in \mathcal{C}^2([0,L])$, the control design objective tackled in this paper is to guarantee the local stability of the system (\ref{eq: wave equation}), when augmented with an adequate control strategy, as well as ensuring the regulation performance, i.e., $z(t) = \frac{\partial y}{\partial x} \rightarrow z_e$ when $t \rightarrow +\infty$. To achieve this objective, we introduce the following deviations: $y_\delta(t,x) = y(t,x) - y_e(x)$ and $u_\delta(t) = u(t) - u_e$. A Taylor expansion with integral remainder shows that (\ref{eq: wave equation}) can equivalently be rewritten under the form:
\begin{subequations}\label{eq: wave equation - variations around equilibrium}
\begin{align}
&\dfrac{\partial^2 y_\delta}{\partial t^2} = \dfrac{\partial^2 y_\delta}{\partial x^2} + f'(y_e) y_\delta + y_\delta^2 \int_0^1 (1-s) f''(y_e + s y_\delta) \,\mathrm{d}s , \label{eq: wave equation - variations around equilibrium - PDE} \\
& y_\delta(t,0)=0,\qquad \dfrac{\partial y_\delta}{ \partial x}(t,L) = u_\delta(t), \label{eq: wave equation - variations around equilibrium - BC} \\
& y_\delta(0,x) = y_0(x)-y_e(x) , \qquad \dfrac{\partial y_\delta}{\partial t}(0,x) = y_1(x) , 
\end{align}
\end{subequations}
for $t > 0$ and $x \in (0,L)$, while the output to be regulated is now expressed as
\begin{equation}\label{eq: system output - variations around equilibrium}
z_\delta(t) = \dfrac{\partial y_\delta}{ \partial x}(t,0) = z(t) - z_e . 
\end{equation}
Finally, following classical proportional integral control design schemes, we introduce the following integral component on the tracking error:
\begin{equation}\label{eq: integral component zeta}
\dot{\zeta}(t) = \dfrac{\partial y_\delta}{ \partial x}(t,0) - z_r(t) = z(t) - ( z_e + z_r(t) ),
\end{equation}
where $z_r(t) \in \R$ is the reference input signal.

\begin{remark}
It was shown in~\cite{lhachemi2019pi} for a linear reaction-diffusion equation with Dirichlet boundary control that a simple proportional integral controller can be used to successfully control a Neumann trace. The control design was performed on a finite-dimensional truncated model capturing the unstable modes of the infinite dimensional system while assessing the stability of the full infinite-dimensional system via a Lyapunov-based argument. Such an approach cannot be directly applied to the case of the wave equation studied in this paper due to the fact that, even in the case of a linear function $f$, the open-loop system might exhibit an infinite number of unstable modes. To avoid this pitfall, we borrow the following remark from~\cite{coron2006global}. In the case  
$f = 0$, the control input $u_\delta(t) = - \alpha \frac{\partial y_{\delta}}{\partial t}(t,L)$, with $\alpha > 0$, ensures the exponential decay of the energy function defined as:
\begin{equation*}
E(t) = \int_0^L \left( \dfrac{\partial y_\delta}{\partial t}(t,x) \right)^2 + \left( \dfrac{\partial y_\delta}{\partial x}(t,x) \right)^2 \,\mathrm{d}x .
\end{equation*}
Thus, as suggested in~\cite{coron2006global}, a suitable control input candidate for (\ref{eq: wave equation - variations around equilibrium}) takes the form:
\begin{equation}\label{eq: preliminary control input}
u_\delta(t) = - \alpha \dfrac{\partial y_\delta}{ \partial t}(t,L) + v(t) ,
\end{equation}
where $\alpha > 0$ is to be selected and $v(t)$ is an auxiliary command input. In particular, it was shown in~\cite{coron2006global} that, in the presence of the nonlinear term $f$, the velocity feedback can be used to locally stabilize all but possibly a finite number of the modes of the system. Then the authors showed that the design of the auxiliary control input $v$ can be performed by pole shifting on a finite dimensional truncated model to achieve the stabilization of the remaining unstable modes. The stability of the resulting closed-loop system was assessed via the introduction of a suitable Lyapunov function. In this paper, we propose to take advantage of such a control design strategy in order to achieve the regulation of the following Neumann trace by means of a proportional integral control design scheme via the introduced integral component (\ref{eq: integral component zeta}).
\end{remark}

\section{Control design}\label{sec: control design}

\subsection{Equivalent homogeneous problem}

By introducing the change of variable:
\begin{equation}\label{eq: change of variable}
w^1(t,x) = y_\delta(t,x) , \qquad 
w^2(t,x) = \dfrac{\partial y_\delta}{\partial t}(t,x) - \dfrac{x}{\alpha L} v(t) ,
\end{equation}
we obtain from the wave equation (\ref{eq: wave equation - variations around equilibrium}), the integral component (\ref{eq: integral component zeta}), and the control strategy (\ref{eq: preliminary control input}) that
\begin{subequations}\label{eq: wave equation - equivalent homogeneous problem}
\begin{align}
& \dfrac{\partial w^1}{\partial t} = w^2 + \dfrac{x}{\alpha L} v(t) ,  \label{eq: wave equation - equivalent homogeneous problem - PDE 1} \\
& \dfrac{\partial w^2}{\partial t} = \dfrac{\partial^2 w^1}{\partial x^2} + f'(y_e) w^1 + r(t,x) - \dfrac{x}{\alpha L} \dot{v}(t) , \label{eq: wave equation - equivalent homogeneous problem - PDE 2} \\
& \dot{\zeta}(t) = \dfrac{\partial w^1}{ \partial x}(t,0) - z_r(t) ,  \\
& w^1(t,0) = 0 ,\qquad \dfrac{\partial w^1}{ \partial x}(t,L) + \alpha w^2(t,L) = 0 , \\
& w^1(0,x) = y_0(x)-y_e(x) , \qquad w^2(0,x) = y_1(x) - \dfrac{x}{\alpha L} v(0) , \\
& \zeta(0) = \zeta_0 
\end{align}
\end{subequations}
for $t > 0$ and $x \in (0,L)$, with the residual term
\begin{equation}\label{eq: def residual term}
r(t,x) = (w^1(t,x))^2 \int_0^1 (1-s) f''(y_e(x) + s w^1(t,x)) \,\mathrm{d}s .
\end{equation}

\begin{remark}
A more classical change of variable for (\ref{eq: wave equation - variations around equilibrium}) with control input $u$ given by (\ref{eq: preliminary control input}) is generally obtained by setting 
\begin{equation}\label{eq: change of variable - classical}
w (t,x) = y_\delta(t,x) - \dfrac{x(x-L)}{L} v(t) .
\end{equation}
In that case, (\ref{eq: wave equation - variations around equilibrium}) with $u$ given by (\ref{eq: preliminary control input}) yields
\begin{subequations}
\begin{align}
& \dfrac{\partial^2 w}{\partial t^2} = \dfrac{\partial^2 w}{\partial x^2} + f'(y_e) w - \dfrac{x(x-L)}{L} \ddot{v}(t) + \left( \dfrac{x(x-L)}{L} f'(y_e) + \dfrac{2}{L} \right) v(t) + r(t,x) , \label{eq: wave equation - equivalent homogeneous problem - classical - PDE} \\
& w(t,0) = 0 , \qquad \dfrac{\partial w}{\partial x}(t,L) + \alpha \dfrac{\partial w}{\partial t}(t,L) = 0 , \\
& w(0,x) = y_0(x) - y_e(x) - \dfrac{x(x-L)}{L} v(0) , \qquad \dfrac{\partial w}{\partial t}(0,x) = y_1(x) - \dfrac{x(x-L)}{L} \dot{v}(0)
\end{align}
\end{subequations}
with 
\begin{align*}
& r(t,x) = \\
& \left( w(t,x) + \dfrac{x(x-L)}{L} v(t) \right)^2 \int_0^1 (1-s) f''\left(y_e(x) + s \left( w(t,x) + \dfrac{x(x-L)}{L} v(t) \right) \right) \,\mathrm{d}s .
\end{align*}
However, this change of variable (\ref{eq: change of variable - classical}) induces the occurrence of a $\ddot{v}$ term in (\ref{eq: wave equation - equivalent homogeneous problem - classical - PDE}), while only a $\dot{v}$ term appears in (\ref{eq: wave equation - equivalent homogeneous problem - PDE 1}-\ref{eq: wave equation - equivalent homogeneous problem - PDE 2}). Thus, in the subsequent procedure, the consideration of the change of variable (\ref{eq: change of variable}) instead of (\ref{eq: change of variable - classical}) will allow a reduction of the complexity of the controller architecture by avoiding the introduction of an extra additional integral component.
\end{remark}

We now introduce the Hilbert space 
\begin{equation*}
\mathcal{H} = \left\{ (w^1,w^2) \in H^1(0,L) \times L^2(0,L) \,:\, w^1(0) = 0 \right\}
\end{equation*}
endowed with the inner product
\begin{equation*}
\left< (w^1,w^2) , (z^1,z^2) \right> 
= \int_0^L (w^1)'(x) \overline{(z^1)'(x)} + w^2(x) \overline{z^2(x)} \,\mathrm{d}x .
\end{equation*}
Defining the following state vector
\begin{equation*}
W(t) = (w^1(t,\cdot),w^2(t,\cdot)) \in \mathcal{H} ,
\end{equation*}
the wave equation with integral component (\ref{eq: wave equation - equivalent homogeneous problem}) can be rewritten under the abstract form
\begin{subequations}\label{eq: wave equation - abstract form}
\begin{align}
& \dfrac{\mathrm{d} W}{\mathrm{d} t}(t) = \mathcal{A} W(t) + a v(t) + b \dot{v}(t) + R(t,\cdot) , \label{eq: wave equation - abstract form - ODE} \\
& \dot{\zeta}(t) = \dfrac{\partial w^1}{ \partial x}(t,0) - z_r(t) , \label{eq: wave equation - abstract form - zeta} \\
& W(0,x) = \left( y_0(x)-y_e(x) , y_1(x) - \dfrac{x}{\alpha L} v(0) \right) , \\
& \zeta(0) = \zeta_0 .
\end{align}
\end{subequations}
for $t > 0$ and $x \in (0,L)$, where 
\begin{equation}\label{eq: def operator A}
\mathcal{A} = 
\begin{pmatrix}
0 & \mathrm{Id} \\ \mathcal{A}_0 & 0
\end{pmatrix}
\end{equation}
with $\mathcal{A}_0 = \Delta + f'(y_e) \,\mathrm{Id}$ on the domain
\begin{align*}
D(\mathcal{A}) 
= \{ (w^1,w^2) \in \mathcal{H} \,:\, & w^1 \in H^2(0,L) ,\, w^2 \in H^1(0,L) ,\, \\
&  w^2(0) = 0 ,\, (w^1)'(L)+\alpha w^2(L) = 0 \} ,
\end{align*} 
and $a,b,R(t,\cdot) \in \mathcal{H}$ are defined by
\begin{equation}\label{eq: def a b R}
\quad a(x) = ( x/(\alpha L) , 0 ) ,
\quad b(x) = ( 0 , -x/(\alpha L) ) , 
\quad R(t,x) = ( 0 , r(t,x) ) .
\end{equation}
We have $R(t,\cdot) \in \mathcal{H}$ because $r(t,\cdot) \in L^2(0,L)$, which is a direct consequence of the facts that $w^1(t,\cdot) \in H^1(0,L) \subset L^\infty(0,L)$, $f''$ is continuous on $\R$, and $y_e$ is continuous on $[0,L]$.

\begin{remark}\label{rmk: classical solutions}
It is well-known that the operator $\mathcal{A}$ generates a $C_0$-semigroup~\cite{tucsnak2009observation}. Moreover, because the Neumann trace is $\mathcal{A}$-admissible, the application of~\cite[Lemma~1]{xu1995robust} shows that the augmentation of $\mathcal{A}$ with the integral component $\zeta$ still generates a $C_0$-semigroup. As $\dot{v}$ is seen as the control input, the state-space vector can further be augmented to include $v$, and the associated augmented operator also generates a $C_0$-semigroup. Now, noting that the residual term (\ref{eq: def residual term}) can be rewritten under the form
\begin{equation*}
r(t,x) = \int_{y_e(x)}^{w^1(t,x)+y_e(x)} (w^1(t,x)+y_e(x)-s) f''(s) \,\mathrm{d}s ,
\end{equation*}
one can see that $w^1 \rightarrow \int_{y_e}^{w^1+y_e} (w^1+y_e-s) f''(s) \,\mathrm{d}s$, when seen as a function from $\{ w^1 \in H^1(0,L) \,:\, w^1(0)=0 \}$ to $L^2(0,L)$, is continuously differentiable. Consequently, the well-posedness of (\ref{eq: wave equation - abstract form}) follows from classical results~\cite{pazy2012semigroups}. In the subsequent developments, we will consider for initial conditions $W(0)\in D(\mathcal{A})$, continuously differentiable reference inputs $z_r$, and a control input $\dot{v}$ that will take the form of a state-feedback, the concept of classical solution for (\ref{eq: wave equation - abstract form}) on its maximal interval of definition $[0,T_{\max})$ with $0 < T_{\max} \leq + \infty$, i.e. $W \in \mathcal{C}^0([0,T_{\max});D(\mathcal{A})) \cap \mathcal{C}^1([0,T_{\max});\mathcal{H})$. 
\end{remark}

\subsection{Properties of the operator $\mathcal{A}$}

First, we explicit in the following lemma the adjoint operator $\mathcal{A}^*$.

\begin{lemma}\label{lem: adjoint operator}
The adjoint operator of $\mathcal{A}$ is defined on
\begin{align*}
D(\mathcal{A}^*) 
= \{ (z^1,z^2) \in \mathcal{H} \,:\, & z^1 \in H^2(0,L) ,\, z^2 \in H^1(0,L) ,\, \\
& z^2(0) = 0 ,\, (z^1)'(L)-\alpha z^2(L) = 0 \} 
\end{align*} 
by
\begin{equation}\label{eq: adjoint operator}
\mathcal{A}^* (z^1,z^2) 
= ( - z^2 - g , - (z^1)'' )
\end{equation}
where $g \in \mathcal{C}^2([0,L])$ is uniquely defined by
\begin{align*}
& g'' = f'(y_e) z^2 , \\
& g(0) = g'(L) = 0 .
\end{align*}
\end{lemma}

\textbf{Proof.}
We write $\mathcal{A} = \mathcal{A}_{tr} + \mathcal{A}_{p}$ with 
\begin{equation*}
\mathcal{A}_{tr} = 
\begin{pmatrix}
0 & \mathrm{Id} \\ \Delta & 0
\end{pmatrix}
,\qquad
\mathcal{A}_{p} = 
\begin{pmatrix}
0 & 0 \\ f'(y_e) \,\mathrm{Id} & 0
\end{pmatrix} 
\end{equation*}
where $\mathcal{A}_{tr}$ is an unbounded operator defined on the same domain as $\mathcal{A}$ while $\mathcal{A}_{p}$ in defined on $\mathcal{H}$. As $\mathcal{A}_{p}$ is bounded, straightforward computations show that $\mathcal{A}_{p}^* = (-g,0)$ where $g$ is defined as in the statement of the Lemma. It remains to compute $\mathcal{A}_{tr}^*$. To do so, one can observe that $0 \in \rho(\mathcal{A}_{tr})$ with
\begin{equation*}
\mathcal{A}_{tr}^{-1} w = \left( - \alpha w^1(L) x + \int_0^x\int_L^\xi w^2(s) \,\mathrm{d}s\,\mathrm{d}\xi , w^1 \right) 
, \quad \forall w = (w^1,w^2) \in \mathcal{H} .
\end{equation*}  
Since $\mathcal{A}_{tr}^{-1}$ is bounded, straightforward computations show that
\begin{equation*}
\left(\mathcal{A}_{tr}^{-1}\right)^* w = \left( - \alpha w^1(L) x - \int_0^x\int_L^\xi w^2(s) \,\mathrm{d}s\,\mathrm{d}\xi , - w^1 \right) 
, \quad \forall w = (w^1,w^2) \in \mathcal{H} .
\end{equation*}
We deduce the claimed result by computing the inverse of the latter operator.
\qed

The strategy reported in this paper relies on the concept of Riesz bases. This concept is recalled in the following definition.

\begin{definition}\label{def: Riesz basis}
A family of vectors $(e_k)_{k \in \Z}$ of $\mathcal{H}$ is a Riesz basis if this family is maximal and there exist constants $m_r,M_R > 0$ such that for any $N \geq 0$ and any $c_{-n_0},\ldots,c_{n_0}\in\C$,
\begin{equation*}
m_R \sum\limits_{\vert k \vert \leq N} \vert c_k \vert^2
\leq \left\Vert \sum\limits_{\vert k \vert \leq N} c_k e_k \right\Vert_\mathcal{H}^2
\leq M_R \sum\limits_{\vert k \vert \leq N} \vert c_k \vert^2 .
\end{equation*}
The dual Riesz basis of $(e_k)_{k \in \Z}$ is the unique family of vectors $(f_k)_{k \in \Z}$ of $\mathcal{H}$ which is such that $\left< e_k , f_l \right>_\mathcal{H} = \delta_{k,l} \in \{0,1\}$ with $\delta_{k,l} = 1$ if and only if $k = l$.
\end{definition}

We can now introduce the following properties of the operator $\mathcal{A}$. These properties, expect the last item, are retrieved from~\cite[Lemmas~2 and~5]{coron2006global}.

\begin{lemma}\label{lem: properties A}
Let $\alpha > 1$. There exists a Riesz basis $(e_k)_{k \in \Z}$ of $\mathcal{H}$ consisting of generalized eigenfunctions of $\mathcal{A}$, associated to the eigenvalues $(\lambda_k)_{k \in \Z}$ and with dual Riesz basis $(f_k)_{k \in \Z}$, such that:
\begin{enumerate}
\item $e_k \in D(\mathcal{A})$ and $\Vert e_k \Vert_\mathcal{H} = 1$ for every $k \in \Z$;
\item each eigenvalue $\lambda_k$ is geometrically simple;
\item there exists $n_0 \geq 0$ such that, for any $k \in \Z$ with $\vert k \vert \geq n_0+1$, the eigenvalue $\lambda_k$  is algebraically simple and satisfies
\begin{equation}\label{eq: lemma properties A - asymptotic behavior eigenvalues}
\lambda_k = \dfrac{1}{2L} \log\left(\dfrac{\alpha-1}{\alpha+1}\right) + i \dfrac{k\pi}{L} + O\left( \dfrac{1}{\vert k \vert} \right)
\end{equation}
as $\vert k \vert \rightarrow + \infty$.
\item if $k \geq n_0+1$, then $e_k$ (resp. $f_k$) is an eigenfunction of $\mathcal{A}$ (resp. $\mathcal{A}^*$) associated with the algebraically simple eigenvalue $\lambda_k$ (resp. $\overline{\lambda_k}$);
\item for every $k \geq n_0+1$, one has $e_k = \overline{e_{-k}}$ and $f_k = \overline{f_{-k}}$;
\item for every $\vert k \vert \leq n_0$, there holds
\begin{subequations}\label{eq: lemma properties A - vectors for small k}
\begin{align}
\mathcal{A} e_k & \in \mathrm{span}\{ e_p \,:\, \vert p \vert \leq n_0 \} , \\
\mathcal{A}^* f_k & \in \mathrm{span}\{ f_p \,:\, \vert p \vert \leq n_0 \} ;
\end{align}
\end{subequations}
\item introducing $e_k = (e_k^1,e_k^2)$, one has $(e_k^1)'(0) = O(1)$ as $\vert k \vert \rightarrow + \infty$.
\end{enumerate}
\end{lemma}

The proof of Lemma~\ref{annex: proof lemma}, which is essentially extracted from~\cite{coron2006global}, is placed in annex for self-completeness of the manuscript.

In the remainder of the paper, we select the constant $\alpha > 1$ such that 
\begin{equation*}
\frac{1}{2L} \log\left(\frac{\alpha-1}{\alpha+1}\right) < -1 .
\end{equation*}
Based on the asymptotic behavior (\ref{eq: lemma properties A - asymptotic behavior eigenvalues}), only a finite number of eigenvalues might have a non negative real part. Thus, without loss of generality, we also select the integer $n_0 \geq 0$ provided by Lemma~\ref{lem: properties A} large enough such that $\operatorname{Re} \lambda_k < -1$ for all $\vert k \vert \geq n_0 + 1$.

\begin{remark}
The state-space can be written as $\mathcal{H} = \mathcal{H}_1 \bigoplus \mathcal{H}_2$ with the subspaces $\mathcal{H}_1 = \mathrm{span}\{ e_p \,:\, \vert p \vert \leq n_0 \}$ and $\mathcal{H}_2 = \overline{\mathrm{span}\{ e_p \,:\, \vert p \vert \geq n_0 + 1 \}}$. Introducing $\pi_1,\pi_2$ the projectors associated with this decomposition, Lemma~\ref{lem: properties A} shows that $\mathcal{A}$ takes the form $\mathcal{A} = \mathcal{A}_1 \pi_1 + \mathcal{A}_2 \pi_2$ where $\mathcal{A}_1 \in \mathcal{L}(\mathcal{H}_1)$ and $\mathcal{A}_2 : D(\mathcal{A}_2) \subset \mathcal{H}_2 \rightarrow \mathcal{H}_2$ with $D(\mathcal{A}_2) = D(\mathcal{A}) \cap \mathcal{H}_2$. Moreover, Lemma~\ref{lem: properties A} shows that $\mathcal{A}_2$ is a Riesz spectral operator. Then, as $D(\mathcal{A}) = \mathcal{H}_1 \bigoplus D(\mathcal{A}_2)$, we obtain that 
\begin{equation*}
D(\mathcal{A})
= \left\{ \sum\limits_{k \in \Z} w_k e_k \,:\, \sum\limits_{k \in \Z} \vert \lambda_k w_k \vert^2 < \infty \right\} 
\end{equation*}
and, for any $w \in D(\mathcal{A})$,
\begin{equation}\label{eq: A = A1 + A2}
\mathcal{A} w = \mathcal{A}_1 \pi_1 w + \sum\limits_{\vert k \vert \geq n_0 + 1} \lambda_k \left< w , f_k \right> e_k ,
\end{equation}
where the equality holds in $\mathcal{H}$-norm. Now, for any $w \in D(\mathcal{A})$, consider the series expansion $w = (w^1,w^2) = \sum_{k \in \Z} \left< w , f_k \right> e_k$. In particular, one has $w^1 = \sum_{k \in \Z} \left< w , f_k \right> e_k^1$ in $H^1$-norm. Moreover, from (\ref{eq: A = A1 + A2}) we have that $\mathcal{A} w = \sum_{k \in \Z} \left< w , f_k \right> \mathcal{A} e_k$ and thus $\mathcal{A}_0 w^1 = \sum_{k \in \Z} \left< w , f_k \right> \mathcal{A}_0 e_k^1$ in $L^2$-norm. Since $f'(y_e) \in L^\infty(0,L)$, the expansion of the latter identity shows that $(w^1)'' = \sum_{k \in \Z} \left< w , f_k \right> (e_k^1)''$ in $L^2$-norm. Consequently, $w^1 = \sum_{k \in \Z} \left< w , f_k \right> e_k^1$ in $H^2$-norm and thus, by the continuous embedding $H^1(0,L) \subset L^\infty(0,L)$,
\begin{equation}\label{eq: series expansion Neumann trace}
(w^1)'(0) = \sum\limits_{k \in \Z} \left< w , f_k \right> (e_k^1)'(0) .
\end{equation}
The latter series expansion will be intensively used in the remainder of this paper. 
\end{remark}

\subsection{Spectral reduction and truncated model}\label{subsec: spectral reduction and truncated model}

Introducing, for every $k \in \Z$,
\begin{equation*}
w_k(t) = \left< W(t) , f_k \right>_\mathcal{H}  , \quad
a_k = \left< a , f_k \right>_\mathcal{H}  , \quad
b_k = \left< b , f_k \right>_\mathcal{H}  , \quad
r_k(t) = \left< R(t,\cdot) , f_k \right>_\mathcal{H} ,
\end{equation*}
we obtain from (\ref{eq: wave equation - abstract form - ODE}), that
\begin{align*}
\dot{w}_k(t)
& = \left< \mathcal{A}W(t) , f_k \right>_\mathcal{H} + a_k v(t) + b_k \dot{v}(t) + r_k(t) \\
& = \left< W(t) , \mathcal{A}^* f_k \right>_\mathcal{H} + a_k v(t) + b_k \dot{v}(t) + r_k(t)
\end{align*}
Recalling that $\mathcal{A}^* f_k = \overline{\lambda_k} f_k$ for $\vert k \vert \geq n_0+1$, we obtain that
\begin{equation}\label{eq: residual dynamics}
\dot{w}_k(t) = \lambda_k w_k(t) + a_k v(t) + b_k \dot{v}(t) + r_k(t) , \qquad \vert k \vert \geq n_0+1 .
\end{equation}
Moreover, after possibly linear recombination\footnote{In that case, all the properties stated by Lemma~\ref{lem: properties A} remain true except that $(e_k)_{\vert k \vert \leq N_0}$ might not be generalized eigenvectors of $\mathcal{A}$.} of $(e_k)_{\vert k \vert \leq N_0}$ and $(f_k)_{\vert k \vert \leq N_0}$, which we still denote by $(e_k)_{\vert k \vert \leq N_0}$ and $(f_k)_{\vert k \vert \leq N_0}$, to obtain matrices with real coefficients, we infer from (\ref{eq: lemma properties A - vectors for small k}) the existence of a matrix $A_0 \in \R^{(2n_0+1)\times(2n_0+1)}$ such that
\begin{equation}\label{eq: prel truncated model}
\dot{X}_0(t) = A_0 X_0(t) + B_{0,1} v(t) + B_{0,2} \dot{v}(t) + R_0(t)
\end{equation}
where
\begin{equation*}
X_0(t) = \begin{bmatrix} w_{-n_0}(t) \\ \vdots \\ w_{n_0}(t) \end{bmatrix} \in \R^{2 n_0 +1} , \quad
B_{0,1} = \begin{bmatrix} a_{-n_0} \\ \vdots \\ a_{n_0} \end{bmatrix} \in \R^{2 n_0 +1} ,
\end{equation*}
\begin{equation*}
B_{0.2} = \begin{bmatrix} b_{-n_0} \\ \vdots \\ b_{n_0} \end{bmatrix} \in \R^{2 n_0 +1} , \quad
R_0(t) = \begin{bmatrix} r_{-n_0}(t) \\ \vdots \\ r_{n_0}(t) \end{bmatrix} \in \R^{2 n_0 +1} .
\end{equation*}
We augment the state-space representation (\ref{eq: prel truncated model}) with the actual control input $v$ as follows:
\begin{equation}\label{eq: truncated model}
\dot{X}_1(t) = A_1 X_1(t) + B_{1} v_d(t) + R_1(t) 
\end{equation}
where
\begin{subequations}
\begin{equation}\label{eq: def X1 vd and R1}
X_1(t) = \begin{bmatrix} v(t) \\ X_{0}(t) \end{bmatrix} \in \R^{2 n_0 +2} , \quad
v_d(t) = \dot{v}(t) , \quad
R_1(t) = \begin{bmatrix} 0 \\ R_0(t) \end{bmatrix} \in \R^{2 n_0 +2}
\end{equation}
\begin{equation}\label{eq: def A1 vd and B1}
A_1 = \begin{bmatrix} 0 & 0 \\ B_{0,1} & A_0 \end{bmatrix} \in \R^{(2n_0+2)\times(2n_0+2)} , \quad
B_1 = \begin{bmatrix} 1 \\ B_{0,2} \end{bmatrix} \in \R^{2 n_0 +2} .
\end{equation}
\end{subequations}

We now further augment the latter state-space representation to include the integral component of the PI controller. First, we note from (\ref{eq: series expansion Neumann trace}) that the dynamics of the integral component $\zeta$ satisfies 
\begin{equation*}
\dot{\zeta}(t) = \sum\limits_{k \in \Z} w_k(t) (e_k^1)'(0) - z_r(t) .
\end{equation*}
We observe that the $\zeta$-dynamics involves all the coefficients of projection $w_k(t)$, with $k \in \Z$, hence cannot be used to augment (\ref{eq: truncated model}). This motivates the introduction of the following change of variable:
\begin{align}
\xi(t) 
& = \zeta(t) - \sum\limits_{\vert k \vert \geq n_0 + 1} \dfrac{(e_k^1)'(0)}{\lambda_k} w_k(t) \label{eq: def integral conponent xi} \\
& = \zeta(t) - 2 \sum\limits_{k \geq n_0 + 1} \operatorname{Re}\left\{ \dfrac{(e_k^1)'(0)}{\lambda_k} w_k(t) \right\} \nonumber
\end{align}
where, using Cauchy-Schwarz inequality, the series is convergent because  $\vert \lambda_k \vert \sim \frac{\vert k \vert \pi}{l}$ and $(e_k^1)'(0) = O(1)$ as $\vert k \vert \rightarrow + \infty$. Moreover, the time derivative of $\xi$ is given by
\begin{align*}
\dot{\xi}(t) 
& = \dot{\zeta}(t) - \sum\limits_{\vert k \vert \geq n_0 + 1} \dfrac{(e_k^1)'(0)}{\lambda_k} \dot{w}_k(t) \\
& = \sum\limits_{\vert k \vert \leq n_0}  w_k(t) (e_k^1)'(0) + \alpha_0 v(t) + \beta_0 \dot{v}(t) - \gamma(t) \\
& = L_1 X_1(t) + \beta_0 v_d(t) - \gamma(t),
\end{align*}
where
\begin{subequations}
\begin{align}
\alpha_0 
& = - \sum\limits_{\vert k \vert \geq n_0 +1} \dfrac{(e_k^1)'(0)}{\lambda_k} a_k
= - 2 \sum\limits_{k \geq n_0 +1} \operatorname{Re}\left\{\dfrac{(e_k^1)'(0)}{\lambda_k} a_k\right\} , \\
\beta_0 
& = - \sum\limits_{\vert k \vert \geq n_0 +1} \dfrac{(e_k^1)'(0)}{\lambda_k} b_k
 = - 2 \sum\limits_{k \geq n_0 +1} \operatorname{Re}\left\{\dfrac{(e_k^1)'(0)}{\lambda_k} b_k\right\} , \\
\gamma(t) 
& = z_r(t) + \sum\limits_{\vert k \vert \geq n_0 +1} \dfrac{(e_k^1)'(0)}{\lambda_k} r_k(t) 
= z_r(t) + 2 \sum\limits_{k \geq n_0 +1} \operatorname{Re}\left\{\dfrac{(e_k^1)'(0)}{\lambda_k} r_k(t)\right\} , \label{eq: def gamma}
\end{align}
\end{subequations}
and
\begin{equation*}
L_1 = \begin{bmatrix} \alpha_0 & (e_{-n_0}^1)'(0) & \ldots & (e_{n_0}^1)'(0) \end{bmatrix} \in \R^{1 \times (2n_0+2)} .
\end{equation*}
Thus, with the introduction of
\begin{subequations}\label{eq: def X A B and Gamma}
\begin{equation}\label{eq: def X and A}
X(t) = \begin{bmatrix} X_1(t) \\ \xi(t) \end{bmatrix} \in \R^{2 n_0 + 3} , \quad
A = \begin{bmatrix} A_1 & 0 \\ L_1 & 0 \end{bmatrix} \in \R^{(2 n_0 + 3)\times(2 n_0 + 3)} ,
\end{equation} 
\begin{equation}\label{eq: def B and Gamma}
B = \begin{bmatrix} B_1 \\ \beta_0 \end{bmatrix} \in \R^{2 n_0 + 3} , \quad
\Gamma(t) = \begin{bmatrix} R_1(t) \\ -\gamma(t) \end{bmatrix} \in \R^{2 n_0 + 3} ,
\end{equation}
\end{subequations}
we obtain the truncated model
\begin{equation}\label{eq: final truncated model}
\dot{X}(t) = A X(t) + B v_d(t) + \Gamma(t) .
\end{equation}

Putting (\ref{eq: residual dynamics}) and (\ref{eq: final truncated model}) together, we obtain that the wave equation with integral component (\ref{eq: wave equation - abstract form}) admits the following equivalent representation used for both control design and stability analysis:
\begin{subequations}\label{eq: wave equation - control design}
\begin{align}
& \dot{X}(t) = A X(t) + B v_d(t) + \Gamma(t) , \label{eq: wave equation - control design - truncated model} \\
& \dot{w}_k(t) = \lambda_k w_k(t) + a_k v(t) + b_k v_d(t) + r_k(t)  , \qquad \vert k \vert \geq n_0+1 . \label{eq: wave equation - residual dynamics}
\end{align}
\end{subequations}

\begin{remark}
The representation (\ref{eq: wave equation - control design}) shows that the dynamics of the wave equation with integral component (\ref{eq: wave equation - abstract form}) can be split into two parts. The first part, given by (\ref{eq: wave equation - control design - truncated model}), consists of an ODE capturing the unstable dynamics plus a certain number of slow stable modes of the system. The second part, given by (\ref{eq: wave equation - residual dynamics}) and referred to as the residual dynamics, captures the stable dynamics of the system which are such that $\operatorname{Re}\lambda_k < -1$. The control strategy consists now into the two following steps. First, a state-feedback is designed to locally stabilize (\ref{eq: wave equation - control design - truncated model}). Then, a stability analysis is carried out to assess that such a control strategy achieves both the local stabilization of (\ref{eq: wave equation - control design}), as well as the fulfillment of the output regulation of the Neumann trace (\ref{eq: system output - variations around equilibrium}).
\end{remark}

\subsection{Control strategy and closed-loop dynamics}\label{sec_poleshift}

The control design strategy consists in the design of a stabilizing state-feedback for (\ref{eq: wave equation - control design - truncated model}). Such a pole shifting is allowed by the following result.

\begin{lemma}\label{lem: Kalman condition}
The pair $(A,B)$ satisfies the Kalman condition.
\end{lemma}

\textbf{Proof.} From (\ref{eq: def X A B and Gamma}), the Hautus test easily shows that $(A,B)$ satisfies the Kalman condition if and only if $(A_1,B_1)$ satisfies the Kalman condition and the square matrix
\begin{equation*}
T = \begin{bmatrix} A_1 & B_1 \\ L_1 & \beta_0 \end{bmatrix} \in \R^{(2 n_0 + 3)\times(2 n_0 + 3)}
\end{equation*}
is invertible. 

We first show the following preliminary result: for any $\lambda \in\C$ and $z \in D(\mathcal{A}^*)$, $\left< a + \lambda b , z \right>_\mathcal{H} = 0$ and $\mathcal{A}^* z = \overline{\lambda} z$ implies $z = 0$. Recall based on Lemma~\ref{lem: adjoint operator} that $\mathcal{A}^* z = \overline{\lambda} z$ gives
\begin{subequations}\label{eq: proof Kalman condition - characterization eigenvector fk}
\begin{align}
& z^2 + g = - \overline{\lambda} z^1 , \label{eq: proof Kalman condition - characterization eigenvector fk - 1} \\
& (z^1)'' = - \overline{\lambda} z^2 , \label{eq: proof Kalman condition - characterization eigenvector fk - 2} \\
& g'' = f'(y_e) z^2 , \label{eq: proof Kalman condition - characterization eigenvector fk - 3} \\
& z^1(0) = z^2(0) = g(0) = g'(L) = 0 , \label{eq: proof Kalman condition - characterization eigenvector fk - 4} \\
& (z^1)'(L)-\alpha z^2(L) = 0 . \label{eq: proof Kalman condition - characterization eigenvector fk - 5}
\end{align}
\end{subequations}
From the definition (\ref{eq: def a b R}) of $a,b\in\mathcal{H}$ one has
\begin{align*}
\left< a + \lambda b , z \right>_\mathcal{H} 
& = \int_0^L \left( \dfrac{x}{\alpha L} \right)' \overline{(z^1)'(x)} - \lambda \dfrac{x}{\alpha L} \overline{z^2(x)} \,\mathrm{dx} \\
& = \left[ \dfrac{x}{\alpha L} \overline{(z^1)'(x)} \right]_{x=0}^{x=L} - \int_0^L \dfrac{x}{\alpha L} \overline{\left( (z^1)''(x) + \overline{\lambda} z^2(x) \right)} \,\mathrm{dx} \\
& = \dfrac{1}{\alpha} \overline{(z^1)'(L)} 
\end{align*} 
where we have used (\ref{eq: proof Kalman condition - characterization eigenvector fk - 2}). Thus we have $(z^1)'(L) = 0$. Then (\ref{eq: proof Kalman condition - characterization eigenvector fk - 5}) shows that $z^2(L) = 0$ and we infer from (\ref{eq: proof Kalman condition - characterization eigenvector fk - 1}) and (\ref{eq: proof Kalman condition - characterization eigenvector fk - 4}) that $(z^2)'(L) = 0$. Moreover, taking twice the derivative of (\ref{eq: proof Kalman condition - characterization eigenvector fk - 1}) and using (\ref{eq: proof Kalman condition - characterization eigenvector fk - 2}-\ref{eq: proof Kalman condition - characterization eigenvector fk - 3}), we obtain that $(z^2)'' + \left( f'(y_e) - (\overline{\lambda})^2 \right) z^2 = 0$. By Cauchy uniqueness, we deduce that $z^2 = 0$. Using (\ref{eq: proof Kalman condition - characterization eigenvector fk - 2}), (\ref{eq: proof Kalman condition - characterization eigenvector fk - 4}) and $(z^1)'(L) = 0$, we reach the conclusion $z = 0$.

Assume now that $(A_1,B_1)$ does not satisfy the Kalman condition. From (\ref{eq: def A1 vd and B1}), the Hautus test shows the existence of $\lambda \in \C$, $x_1 \in \C$, and $x_2 \in \C^{2n_0+1}$, with either $x_1 \neq 0$ or $x_2 \neq 0$, such that
\begin{align*}
& x_2^* B_{0,1} = \lambda x_1^* , \\
& x_2^* A_0 = \lambda x_2^* , \\
& x_1^* + x_2^* B_{0,2} = 0 .
\end{align*}
This implies the existence of $x_2 \neq 0$ such that
\begin{align*}
& A_0^* x_2 = \overline{\lambda} x_2 , \\
& x_2^* ( B_{0,1} + \lambda B_{0,2} ) = 0 
\end{align*}
where $B_{0,1} + \lambda B_{0,2} = ( \left< a + \lambda b , f_k \right>_\mathcal{H} )_{-n_0 \leq k \leq n_0}$. Noting that $A_0^*$ is the matrix of $\mathcal{A}^*$ in $(f_k)_{\vert k \vert \leq n_0}$, this shows the existence of a nonzero vector $z \in D(\mathcal{A}^*)$ such that $\left< a + \lambda b , z \right>_\mathcal{H} = 0$ and $\mathcal{A}^* z = \overline{\lambda} z$. The result of the previous paragraph leads to the contraction $z=0$. Hence $(A_1,B_1)$ does satisfy the Kalman condition. 

It remains to show that the matrix $T$ is invertible. Let a vector 
\begin{equation*}
X_e = \begin{bmatrix}
v_e & w_{-n_0,e} & \ldots & w_{-n_0,e} & v_{d,e} 
\end{bmatrix}^\top \in \R^{2n_0+3}
\end{equation*}
be an element of the kernel of $T$. Then, by expanding $T X_e = 0$, we obtain that
\begin{align*}
0 & = v_{d,e} , \\
0 & = A_0 \begin{bmatrix} w_{-n_0,e} \\ \vdots \\ w_{-n_0,e} \end{bmatrix} + B_{0,1} v_e , \\
0 & =  \sum\limits_{\vert k \vert \leq n_0} w_{k,e} (e_k^1)'(0) - \left( \sum\limits_{\vert k \vert \geq n_0 +1} \dfrac{(e_k^1)'(0)}{\lambda_k} a_k \right) v_e .
\end{align*}
We define, for $\vert k \vert \geq n_0+1$, $w_{k,e} = - \frac{a_k}{\lambda_k} v_e$. Then, as $(w_{k,e})_{k \in \Z}$ is square summable, we can introduce $w_e = \sum_{k \in \Z} w_{k,e} e_k \in \mathcal{H}$. We obtain that
\begin{align*}
0 & = A_0 \begin{bmatrix} w_{-n_0,e} \\ \vdots \\ w_{-n_0,e} \end{bmatrix} + B_{0,1} v_e , \\
0 & = \lambda_k w_{k,e} + a_k v_e , \quad \vert k \vert \geq n_0+1 , \\
0 & =  \sum\limits_{k \in \Z} w_{k,e} (e_k^1)'(0) .
\end{align*}
In particular $(\lambda_k w_{k,e})_{k \in \Z}$ is square summable and thus $w_e \in D(\mathcal{A})$. The developments of Subsection~\ref{subsec: spectral reduction and truncated model} show that the above system is equivalent to 
\begin{align*}
& \mathcal{A}w_e + a v_e = 0 , \\
& (w_e^1)'(0) = 0 
\end{align*}
with $w_e = (w_e^1,w_e^2)$. By expanding the former identity, we first have that $(w_e^1)''+f'(y_e)w_e^1=0$ with $w_e^1(0)=(w_e^1)'(0)=0$ and thus, by Cauchy uniqueness, $w_e^1 = 0$. We also have $w_e^2 = \frac{-x}{\alpha L} v_e$ with $v_e = - \alpha w_e^2(L) = (w_e^1)'(L) = 0$, and thus $w_e^2 = 0$. This yields $w_e = 0$, which shows that $w_{k,e} = 0$ for every $k \in \Z$. Hence the kernel of $T$ is reduced to $\{0\}$, which concludes the proof. 
\qed

The result of Lemma~\ref{lem: Kalman condition} ensures the existence of a gain $K \in \R^{1 \times (2 n_0 +3)}$ such that $A_K = A + BK$ is Hurwitz. Then, we can set the state feedback
\begin{equation}\label{eq: control input vd}
v_d(t) = \dot{v}(t) = K X(t) ,
\end{equation}
which yields the following closed-loop system dynamics:
\begin{subequations}\label{eq: wave equation - closed-loop}
\begin{align}
& \dot{X}(t) = A_K X(t) + \Gamma(t) ,  \\
& \dot{w}_k(t) = \lambda_k w_k(t) + a_k v(t) + b_k v_d(t) + r_k(t)  , \qquad \vert k \vert \geq n_0+1 .
\end{align}
\end{subequations}
The objective of the remainder of the paper is to assess the local stability of the closed-loop system, as well as the study of the tracking performance.

\section{Stability and set point regulation assessment}\label{sec: stability and regulation}

\subsection{Stability analysis}

The main stability result of this paper is stated in the following theorem.

\begin{theorem} \label{thm: stability result}
There exist $\kappa \in (0,1)$ and $\overline{C}_1,\delta > 0$ such that, for any $\eta \in [0,1)$, there exists $\overline{C}_2 > 0$ such that, for any initial condition satisfying
\begin{equation*}
\Vert W(0) \Vert_\mathcal{H}^2 + \vert \xi(0) \vert^2 + \vert v(0) \vert^2 \leq \delta
\end{equation*}
and any continuously differentiable reference input $z_r$ with
\begin{equation*}
\Vert z_r \Vert_{L^\infty(\R_+)}^2 \leq \delta , 
\end{equation*}
the classical solutions of (\ref{eq: wave equation - equivalent homogeneous problem}) with control law (\ref{eq: control input vd}) is well defined on $\R_+$ and satisfies
\begin{equation} \label{eq: stability result - w^1 L^inf bounded by one}
\Vert w^1(t,\cdot) \Vert_{L^\infty(0,L)} < 1
\end{equation}
and 
\begin{align}
\Vert W(t) \Vert_\mathcal{H}^2 + \vert \xi(t) \vert^2 + \vert v(t) \vert^2 
& \leq \overline{C}_1 e^{- 2 \kappa t} \left( \Vert W(0) \Vert_\mathcal{H}^2 + \vert \xi(0) \vert^2 + \vert v(0) \vert^2 \right) \label{eq: stability result - local exp stab} \\
& \phantom{\leq}\, + \overline{C}_2 \sup\limits_{0 \leq s \leq t} e^{-2\eta\kappa(t-s)} \vert z_r(s) \vert^2 . \nonumber
\end{align}
for all $t \geq 0$.
\end{theorem}

\begin{remark}
In the particular case of a linear function $f$, which implies that the residual term (\ref{eq: def residual term}) is identically zero, the exponential stability result (\ref{eq: stability result - local exp stab}) stated by Theorem~\ref{thm: stability result} is global.
\end{remark}

\begin{remark} The result of Theorem~\ref{thm: stability result} ensures the stability of the closed-loop system in $(w,\xi)$ coordinates. This immediately induces the stability of the closed-loop system in its original coordinates because, from (\ref{eq: change of variable}), we have
\begin{align*}
\left\Vert \left( y_\delta(t,\cdot) , \dfrac{\partial y_\delta}{\partial t}(t,\cdot) \right) \right\Vert_\mathcal{H}
& \leq \Vert W(t) \Vert_\mathcal{H} + \left\Vert \left( 0 , \dfrac{(\cdot)}{\alpha L} v(t) \right) \right\Vert_\mathcal{H} \\
& \leq \Vert W(t) \Vert_\mathcal{H} + \dfrac{1}{\alpha} \sqrt{\dfrac{L}{3}} \vert v(t) \vert 
\end{align*}
and, from (\ref{eq: def integral conponent xi}),
\begin{align*}
\vert \zeta(t) \vert 
& \leq \vert \xi(t) \vert + \sum\limits_{\vert k \vert \geq n_0 + 1} \left\vert \dfrac{(e_k^1)'(0)}{\lambda_k} w_k(t) \right\vert \\
& \leq \vert \xi(t) \vert + \sqrt{\sum\limits_{\vert k \vert \geq n_0 + 1} \left\vert \dfrac{(e_k^1)'(0)}{\lambda_k} \right\vert^2} \sqrt{\sum\limits_{\vert k \vert \geq n_0 + 1} \vert w_k(t) \vert^2} \\
& \leq \vert \xi(t) \vert + \sqrt{\dfrac{1}{m_R} \sum\limits_{\vert k \vert \geq n_0 + 1} \left\vert \dfrac{(e_k^1)'(0)}{\lambda_k} \right\vert^2} \Vert W(t) \Vert_\mathcal{H}   
\end{align*}
where we recall that $(e_k^1)'(0) = O(1)$ and $\vert \lambda_k \vert \sim k\pi/L$ when $\vert k \vert \rightarrow + \infty$.
\end{remark}

\textbf{Proof of Theorem~\ref{thm: stability result}.} Let $M > 3( \Vert a \Vert_\mathcal{H}^2 + \Vert b \Vert_\mathcal{H}^2 \Vert K \Vert^2 )/m_R$ be given, where we recall that $a,b \in \mathcal{H}$ are defined by (\ref{eq: def a b R}) and the constant $m_R > 0$ is as provided by Definition~\ref{def: Riesz basis}. We introduce for all $t \geq 0$
\begin{equation}\label{eq: def Lyapunov func}
V(t) = M X(t)^\top P X(t) + \dfrac{1}{2} \sum\limits_{\vert k \vert \geq n_0+1} \vert w_k(t) \vert^2 ,
\end{equation}
where $P$ is a symmetric definite positive matrix such that $A_K^\top P + P A_K = -I$. Then we obtain from (\ref{eq: wave equation - closed-loop}) that
\begin{align*}
\dot{V}(t)
& = M X(t)^\top \left\{ A_K^\top P + P A_K \right\} X(t) + M \left\{ \Gamma(t)^\top P X(t) + X(t)^\top P \Gamma(t) \right\} \\
& \phantom{=}\, + \sum\limits_{\vert k \vert \geq n_0+1} \operatorname{Re}\lambda_k \vert w_k(t) \vert^2 + \sum\limits_{\vert k \vert \geq n_0+1} \operatorname{Re} \left\{ \overline{w_k(t)} \left( a_k v(t) + b_k v_d(t) + r_k(t) \right) \right\} \\
& \leq - M \Vert X(t) \Vert^2 - \sum\limits_{\vert k \vert \geq n_0+1} \vert w_k(t) \vert^2 + 2 M \Vert P \Vert \Vert X(t) \Vert \Vert \Gamma(t) \Vert  \\
& \phantom{=}\, + \sum\limits_{\vert k \vert \geq n_0+1} \vert w_k(t) \vert \left( \vert a_k \vert \vert v(t) \vert + \vert b_k \vert \vert v_d(t) \vert + \vert r_k(t) \vert \right)
\end{align*}
where we have used that $\operatorname{Re} \lambda_k < -1$ for all $\vert k \vert \geq n_0 + 1$. Using now Young's inequality, we infer that
\begin{align*}
2 M \Vert P \Vert \Vert X(t) \Vert \Vert \Gamma(t) \Vert
& \leq \dfrac{M}{2} \Vert X(t) \Vert ^2 + 2 M \Vert P \Vert^2 \Vert \Gamma(t) \Vert^2
\end{align*}
and
\begin{align*}
& \sum\limits_{\vert k \vert \geq n_0+1} \vert w_k(t) \vert \left( \vert a_k \vert \vert v(t) \vert + \vert b_k \vert \vert v_d(t) \vert + \vert r_k(t) \vert \right) \\
& \qquad \leq \dfrac{1}{2} \sum\limits_{\vert k \vert \geq n_0+1} \vert w_k(t) \vert^2 + \dfrac{3}{2} \sum\limits_{\vert k \vert \geq n_0+1}  \left( \vert a_k \vert^2 \vert v(t) \vert^2 + \vert b_k \vert^2 \vert v_d(t) \vert^2 + \vert r_k(t) \vert^2 \right) \\
& \qquad \leq \dfrac{1}{2} \sum\limits_{\vert k \vert \geq n_0+1} \vert w_k(t) \vert^2 + \dfrac{3\Vert a \Vert_\mathcal{H}^2}{2m_R} \vert v(t) \vert^2 + \dfrac{3\Vert b \Vert_\mathcal{H}^2}{2m_R} \vert v_d(t) \vert^2 + \dfrac{3}{2m_R} \Vert R(t,\cdot) \Vert_\mathcal{H}^2 \\
& \qquad \leq \dfrac{1}{2} \sum\limits_{\vert k \vert \geq n_0+1} \vert w_k(t) \vert^2 + \dfrac{3 (\Vert a \Vert_\mathcal{H}^2 + \Vert b \Vert_\mathcal{H}^2 \Vert K \Vert^2)}{2m_R} \Vert X(t) \Vert^2 + \dfrac{3}{2m_R} \Vert R(t,\cdot) \Vert_\mathcal{H}^2 
\end{align*}
where we have used (\ref{eq: control input vd}) and the fact that $v(t)$ is the first component of $X(t)$. Thus, we obtain that 
\begin{align}
\dot{V}(t)
& \leq - \left\{ \dfrac{M}{2} - \dfrac{3 (\Vert a \Vert_\mathcal{H}^2 + \Vert b \Vert_\mathcal{H}^2 \Vert K \Vert^2)}{2m_R} \right\} \Vert X(t) \Vert^2 - \dfrac{1}{2} \sum\limits_{\vert k \vert \geq n_0+1} \vert w_k(t) \vert^2 \label{eq: prel estimate dot_V} \\
& \phantom{=}\, + 2 M \Vert P \Vert^2 \Vert \Gamma(t) \Vert^2 + \dfrac{3}{2m_R} \Vert R(t,\cdot) \Vert_\mathcal{H}^2 . \nonumber
\end{align}
We evaluate the two last terms of the above inequality. Recalling that $R_1(t)$, $\gamma(t)$, and $\Gamma(t)$ are defined by (\ref{eq: def X1 vd and R1}), (\ref{eq: def gamma}), and (\ref{eq: def B and Gamma}), respectively, we have
\begin{align}
\Vert \Gamma(t) \Vert^2 
& = \Vert R_1(t) \Vert^2 + \vert \gamma(t) \vert^2 \nonumber \\
& = \sum\limits_{\vert k \vert \leq n_0} \vert r_k(t) \vert^2 + \left\vert z_r(t) + \sum\limits_{\vert k \vert \geq n_0 +1} \dfrac{(e_k^1)'(0)}{\lambda_k}r_k(t) \right\vert^2 \nonumber \\
& \leq \sum\limits_{\vert k \vert \leq n_0} \vert r_k(t) \vert^2 + 2 \vert z_r(t) \vert^2 + 2 \sum\limits_{\vert k \vert \geq n_0 +1} \left\vert\dfrac{(e_k^1)'(0)}{\lambda_k}\right\vert^2  \sum\limits_{\vert k \vert \geq n_0 +1} \vert r_k(t) \vert^2 \nonumber \\
& \leq C_0^2 \Vert R(t,\cdot) \Vert_\mathcal{H}^2 + 2 \vert z_r(t) \vert^2 \label{eq: estimate Gamma}
\end{align}
with $C_0 > 0$ defined by 
\begin{equation*}
C_0^2 = \frac{1}{m_R} \max\left( 1 , 2 \sum\limits_{\vert k \vert \geq n_0 +1} \left\vert\dfrac{(e_k^1)'(0)}{\lambda_k}\right\vert^2 \right) < + \infty .
\end{equation*}
We now evaluate $\Vert R(t,\cdot) \Vert_\mathcal{H}^2 = \int_0^L \vert r(t,x) \vert^2 \,\mathrm{d}x $ where we recall that $r(t,x)$ is defined by (\ref{eq: def residual term}). Let $\epsilon \in (0,1)$ to be determined. Let $C_I \geq 0$ be the maximum of $\vert f'' \vert$ over the range $\left[ \min y_e -1 , \max y_e +1 \right]$. Thus, assuming that 
\begin{equation}\label{eq: a priori estimate}
\Vert w^1(t,\cdot) \Vert_{L^\infty(0,L)} \leq \epsilon < 1 ,
\end{equation}
we obtain that $\vert r(t,x) \vert \leq C_I \vert w^1(t,x) \vert^2$ and thus
\begin{align}
\Vert R(t,\cdot) \Vert_\mathcal{H}^2
& = \int_0^L \vert r(t,x) \vert^2 \,\mathrm{d}x \nonumber \\
& \leq C_I^2 \int_0^L \vert w^1(t,x) \vert^4 \,\mathrm{d}x \nonumber \\
& \leq \epsilon^2 C_I^2 \Vert w^1(t,\cdot) \Vert_{L^2(0,L)}^2 \nonumber \\
& \leq \epsilon^2 L^2 C_I^2 \Vert W(t) \Vert_{\mathcal{H}}^2 \label{eq: prel estimate R} \\
& \leq \epsilon^2 M_R L^2 C_I^2 \sum\limits_{k \in \Z} \vert w_k(t) \vert^2 \nonumber \\
& \leq \epsilon^2 C_1^2 \left\{ \Vert X(t) \Vert^2 + \sum\limits_{\vert k \vert \geq n_0+1} \vert w_k(t) \vert^2 \right\} \label{eq: estimate R}
\end{align}
where the third inequality follows from Poincar{\'e} inequality, the fourth inequality follows from Definition~\ref{def: Riesz basis}, and with the constant $C_1 \geq 0$ defined by $C_1^2 = M_R L^2 C_I^2$. We deduce from (\ref{eq: prel estimate dot_V}-\ref{eq: estimate Gamma}) and, under the \emph{a priori} estimate (\ref{eq: a priori estimate}), (\ref{eq: estimate R}) that
\begin{align*}
\dot{V}(t)
& \leq - \left\{ \dfrac{M}{2} - \dfrac{3 (\Vert a \Vert_\mathcal{H}^2 + \Vert b \Vert_\mathcal{H}^2 \Vert K \Vert^2)}{2m_R} - \epsilon^2 C_1^2 \left( 2M\Vert P \Vert^2 C_0^2 + \dfrac{3}{2 m_R} \right) \right\} \Vert X(t) \Vert^2 \\
& \phantom{=}\, - \dfrac{1}{2} \left\{( 1 - 2\epsilon^2 C_1^2 \left( 2M\Vert P \Vert^2 C_0^2 + \dfrac{3}{2 m_R} \right) \right\} \sum\limits_{\vert k \vert \geq n_0+1} \vert w_k(t) \vert^2  \\
& \phantom{=}\, + 4 M \Vert P \Vert^2 \vert z_r(t) \vert^2 .
\end{align*}
With $M > 3( \Vert a \Vert_\mathcal{H}^2 + \Vert b \Vert_\mathcal{H}^2 \Vert K \Vert^2 )/m_R$ and by selecting $\epsilon \in (0,1)$ small enough (independently of the initial conditions and the reference signal $z_r$), we obtain the existence of constants $\kappa,C_2 > 0$ such that, under the \emph{a priori} estimate (\ref{eq: a priori estimate}),
\begin{align*}
\dot{V}(t)
& \leq - 2 \kappa V(t) + C_2 \vert z_r(t) \vert^2 .
\end{align*}
Let $\eta\in[0,1)$ be arbitrary. Assuming that the \emph{a priori} estimate (\ref{eq: a priori estimate}) holds on $[0,t]$ for some $t > 0$, we obtain that
\begin{align*}
V(t)
& \leq e^{- 2 \kappa t} V(0) + C_2 \int_0^t e^{-2\kappa(t-s)} \vert z_r(s) \vert^2 \,\mathrm{d}s \\
& \leq e^{- 2 \kappa t} V(0) + \dfrac{C_2}{2\kappa(1-\eta)} \sup\limits_{0 \leq s \leq t} e^{-2\eta\kappa(t-s)} \vert z_r(s) \vert^2 .
\end{align*}
Denoting by $\lambda_m(P) > 0$ and $\lambda_M(P) > 0$ the smallest and largest eigenvalues of the symmetric definite positive matrix $P$, we now note from (\ref{eq: def Lyapunov func}) that 
\begin{align*}
V(0) & 
\leq M \lambda_M(P) \Vert X(0) \Vert^2 + \dfrac{1}{2} \sum\limits_{\vert k \vert \geq n_0+1} \vert w_k(0) \vert^2 \\
& \leq C_3 \left( \Vert W(0) \Vert_\mathcal{H}^2 + \vert \xi(0) \vert^2 + \vert v(0) \vert^2 \right)
\end{align*}
for some constant $C_3 > 0$ and 
\begin{align*}
V(t) & 
\geq M \lambda_m(P) \Vert X(t) \Vert^2 + \dfrac{1}{2} \sum\limits_{\vert k \vert \geq n_0+1} \vert w_k(t) \vert^2 \\
& \geq C_4 \left( \Vert W(t) \Vert_\mathcal{H}^2 + \vert \xi(t) \vert^2 + \vert v(t) \vert^2 \right)
\end{align*}
for some constant $C_4 > 0$. Thus, assuming that the \emph{a priori} estimate (\ref{eq: a priori estimate}) holds over $[0,t]$ for some $t > 0$,  we deduce from the three above estimates the existence of constants $\overline{C}_1,\overline{C}_2 > 0$ such that 
\begin{align}
\Vert W(t) \Vert_\mathcal{H}^2 + \vert \xi(t) \vert^2 + \vert v(t) \vert^2 
& \leq \overline{C}_1 e^{- 2 \kappa t} \left( \Vert W(0) \Vert_\mathcal{H}^2 + \vert \xi(0) \vert^2 + \vert v(0) \vert^2 \right) \label{eq: stab result} \\
& \phantom{\leq}\, + \overline{C}_2 \sup\limits_{0 \leq s \leq t} e^{-2\eta\kappa(t-s)} \vert z_r(s) \vert^2 . \nonumber
\end{align}
To conclude, let us note that 
\begin{equation*}
\Vert w^1(t,\cdot) \Vert_{L^\infty(0,L)}
\leq \sqrt{L} \Vert (w^1)'(t,\cdot) \Vert_{L^2(0,L)}
\leq \sqrt{L} \Vert W(t) \Vert_\mathcal{H} .
\end{equation*}
Hence, if the initial condition is selected such that 
\begin{equation*}
\Vert W(0) \Vert_\mathcal{H}^2 + \vert \xi(0) \vert^2 + \vert v(0) \vert^2 \leq \dfrac{\epsilon^2}{2L} \min\left( 1 , \frac{1}{2 \overline{C}_1} \right) ,
\end{equation*}
which in particular ensures that $\Vert w^1(0,\cdot) \Vert_{L^\infty(0,L)} \leq \epsilon/\sqrt{2} < \epsilon$, and the reference input is chosen such that 
\begin{equation*}
\Vert z_r \Vert_{L^\infty(0,L)}^2 \leq \frac{\epsilon^2}{4 L\overline{C}_2} ,
\end{equation*}
it is readily checked based on (\ref{eq: stab result}) that the \emph{a priori} estimate (\ref{eq: a priori estimate}) holds for all $t \geq 0$. In this case, the stability estimate (\ref{eq: stab result}) holds for all $t \geq 0$.
\qed

\subsection{Set point regulation}

We are now in position to assess the set point regulation of the left Neumann trace.

\begin{theorem}\label{thm: reference tracking}
Let $\kappa \in (0,1)$, $\delta > 0$, and $\eta \in [0,1)$ be as provided by Theorem~\ref{thm: stability result}. There exist constants $\overline{C}_3,\overline{C}_4>0$ such that, for any initial condition satisfying
\begin{equation*}
\Vert W(0) \Vert_\mathcal{H}^2 + \vert \xi(0) \vert^2 + \vert v(0) \vert^2 \leq \delta
\end{equation*}
and any continuously differentiable reference input $z_r$ with
\begin{equation*}
\Vert z_r \Vert_{L^\infty(\R_+)}^2 \leq \delta , 
\end{equation*}
the classical solutions of (\ref{eq: wave equation - equivalent homogeneous problem}) satisfies
\begin{align}
\left\vert (w^1)'(t,0) - z_r(t) \right\vert
& \leq \overline{C}_3 e^{-\kappa t} \left\{ \Vert W(0) \Vert_\mathcal{H} + \vert \xi(0) \vert + \vert v(0) \vert + \Vert \mathcal{A} W(0) \Vert_\mathcal{H} \right\} \label{eq: thm ref tracking - estimate} \\
& \phantom{\leq}\; + \overline{C}_4 \sup\limits_{0 \leq s \leq t} e^{-\eta\kappa(t-s)} \vert z_r(s) \vert \nonumber
\end{align}
for all $t \geq 0$.
\end{theorem}

\begin{remark}
In particular, $z_r(t) \rightarrow 0$ implies $(w^1)'(t,0) \rightarrow 0$, i.e. $z(t) \rightarrow z_e$, which achieves the desired set point reference tracking.
\end{remark}

Before proceeding with the proof of Theorem~\ref{thm: reference tracking}, we first derive an estimate of $\Vert \frac{\mathrm{d}R}{\mathrm{d}t}(t,\cdot) \Vert_\mathcal{H}$. To do so, we assume that the assumptions and conclusions of Theorem~\ref{thm: stability result} apply. Following up with Remark~\ref{rmk: classical solutions}, we have that $\frac{\mathrm{d}R}{\mathrm{d}t} = (0,\frac{\mathrm{d}r}{\mathrm{d}t})$ with
\begin{equation}\label{eq: time derivative r}
\dfrac{\mathrm{d}r}{\mathrm{d}t}(t,\cdot) = \left[ w^2(t,\cdot) + \dfrac{(\cdot)}{\alpha L} v(t) \right] \int_{y_e}^{w^1(t,\cdot)+y_e} f''(s) \,\mathrm{d}s .
\end{equation}
Using (\ref{eq: stability result - w^1 L^inf bounded by one}) and the estimate $\vert f'' \vert \leq C_I$ on the range $\left[ \min y_e -1 , \max y_e +1 \right]$, we deduce that 
\begin{align}
\left\Vert \dfrac{\mathrm{d}R}{\mathrm{d}t}(t,\cdot) \right\Vert_\mathcal{H}^2
& = \left\Vert \dfrac{\mathrm{d}r}{\mathrm{d}t}(t,\cdot) \right\Vert_{L^2(0,L)}^2 \nonumber \\
& \leq \int_0^L \left\vert w^2(t,x) + \dfrac{x}{\alpha L} v(t) \right\vert^2 \left\vert \int_{y_e(x)}^{w^1(t,x)+y_e(x)} \vert f''(s) \vert \,\mathrm{d}s \right\vert^2 \,\mathrm{d}x \nonumber \\
& \leq 2 C_I^2 \left\{ \Vert w^2(t) \Vert_{L^2(0,L)}^2 + \dfrac{L}{3 \alpha^2} \vert v(t) \vert^2 \right\} \nonumber \\
& \leq 2 C_I^2 \left\{ \Vert W(t) \Vert_\mathcal{H}^2 + \dfrac{L}{3 \alpha^2} \vert v(t) \vert^2 \right\}. \label{eq: estimate norm R_t}
\end{align}
We are now in position to complete the proof of Theorem~\ref{thm: reference tracking}.

\textbf{Proof of Theorem~\ref{thm: reference tracking}.} We fix an integer $N \geq n_0$ and a constant $\gamma > 0$ such that\footnote{We recall that $\operatorname{Re}\lambda_k < -1$ for all $\vert k \vert \geq n_0 +1$ and that $\kappa \in (0,1)$.} $\mathrm{Re}\lambda_k \leq -\gamma < -\kappa < 0$ for all $\vert k \vert \geq N+1$. Then we have from (\ref{eq: series expansion Neumann trace}) that
\begin{align}
& \left\vert (w^1)'(t,0) - z_r(t) \right\vert \nonumber \\
& \qquad\leq \sum\limits_{k \in \Z} \vert w_k(t) \vert \vert (e_k^1)'(0) \vert + \vert z_r(t) \vert \nonumber \\
& \qquad\leq \sum\limits_{\vert k \vert \leq N} \vert w_k(t) \vert \vert (e_k^1)'(0) \vert + \sum\limits_{\vert k \vert \geq N+1} \vert \lambda_k w_k(t) \vert \left\vert \dfrac{(e_k^1)'(0)}{\lambda_k} \right\vert + \vert z_r(t) \vert \nonumber \\
& \qquad\leq \sqrt{\sum\limits_{\vert k \vert \leq N} \vert (e_k^1)'(0) \vert^2} \sqrt{\sum\limits_{\vert k \vert \leq N} \vert w_k(t) \vert^2} \nonumber \\
& \qquad\phantom{\leq}\; + \sqrt{\sum\limits_{\vert k \vert \geq N+1} \left\vert \dfrac{(e_k^1)'(0)}{\lambda_k} \right\vert^2} \sqrt{\sum\limits_{\vert k \vert \geq N+1} \vert \lambda_k w_k(t) \vert^2} + \vert z_r(t) \vert \nonumber \\
& \qquad\leq \sqrt{\dfrac{1}{m_R} \sum\limits_{\vert k \vert \leq N} \vert (e_k^1)'(0) \vert^2} \Vert W(t) \Vert_\mathcal{H} \label{eq: prel estimate tracking performance} \\
& \qquad\phantom{\leq}\; + \sqrt{\sum\limits_{\vert k \vert \geq N+1} \left\vert \dfrac{(e_k^1)'(0)}{\lambda_k} \right\vert^2} \sqrt{\sum\limits_{\vert k \vert \geq N+1} \vert \lambda_k w_k(t) \vert^2} + \vert z_r(t) \vert \nonumber
\end{align}
Based on (\ref{eq: stability result - local exp stab}), we only need to evaluate the term $\sum_{\vert k \vert \leq N} \vert \lambda_k w_k(t) \vert^2$. We have from (\ref{eq: wave equation - closed-loop}) that, for all $\vert k \vert \geq N+1 \geq n_0+1$,
\begin{align}
\lambda_k w_k(t) 
& = e^{\lambda_k t} \lambda_k w_k(0) + \int_0^t \lambda_k e^{\lambda_k(t-s)} \left\{ a_k v(s) + b_k v_d(s) + r_k(s) \right\} \,\mathrm{d}s \nonumber \\
& = e^{\lambda_k t} \lambda_k w_k(0) - \left\{ a_k v(t) + b_k v_d(t) + r_k(t) \right\} \label{eq: lambda_k w_k}\\
& \phantom{=}\, + e^{\lambda_k t} \left\{ a_k v(0) + b_k v_d(0) + r_k(0) \right\} \nonumber \\
& \phantom{=}\, + \int_0^t e^{\lambda_k(t-s)} \left\{ a_k v_d(s) + b_k \dot{v}_d(s) + \dot{r}_k(s) \right\} \,\mathrm{d}s \nonumber .
\end{align}
We have from (\ref{eq: control input vd}-\ref{eq: wave equation - closed-loop}) that $v_d(t) = K X(t)$ and $\dot{v}_d(t) = K A_K X(t) + K \Gamma(t)$ with 
\begin{subequations}\label{eq: set point - intermediate estimates for lamda_k wk}
\begin{align}
\Vert X(t) \Vert^2 & \leq \dfrac{1}{m_R} \Vert W(t) \Vert_\mathcal{H}^2 + \vert \xi(t) \vert^2 + \vert v(t) \vert^2 , \\
\Vert \Gamma(t) \Vert^2 & \leq L^2 C_0^2 C_I^2 \Vert W(t) \Vert_\mathcal{H}^2 + 2 \vert z_r(t) \vert^2 
\end{align}
where the second inequality follows from (\ref{eq: estimate Gamma}) and (\ref{eq: prel estimate R}). For ease of notation, we define $\mathrm{CI} = \sqrt{\Vert W(0) \Vert_\mathcal{H}^2 + \vert \xi(0) \vert^2 + \vert v(0) \vert^2}$. Using $\gamma > \kappa > \eta\kappa \geq 0$, we have  from (\ref{eq: stability result - local exp stab}) that
\begin{align}
& \left\vert \int_0^t e^{\lambda_k(t-s)} v_d(s) \,\mathrm{d}s \right\vert \nonumber \\
& \leq \Vert K \Vert \int_0^t e^{-\gamma(t-s)} \Vert X(s) \Vert \,\mathrm{d}s \nonumber \\
& \leq \Vert K \Vert \max(1,m_R^{-1/2}) e^{-\gamma t} \int_0^t e^{\gamma s} \left\{ \sqrt{\overline{C}_1} e^{-\kappa s}\mathrm{CI} + \sqrt{\overline{C}_2} \sup\limits_{0 \leq \tau \leq s} e^{-\eta\kappa(s-\tau)} \vert z_r(\tau) \vert \right\} \,\mathrm{d}s \nonumber\\
& \leq C_5 e^{-\kappa t} \mathrm{CI} + C_5 \sup\limits_{0 \leq \tau \leq t} e^{-\eta\kappa(t-\tau)} \vert z_r(\tau) \vert \label{eq: set point - intermediate estimates for lamda_k wk - 1}
\end{align}
for some constant $C_5 > 0$ and, similarly,
\begin{align}
\left\vert \int_0^t e^{\lambda_k(t-s)} \dot{v}_d(s) \,\mathrm{d}s \right\vert
& \leq C_6 e^{-\kappa t} \mathrm{CI} + C_6 \sup\limits_{0 \leq \tau \leq t} e^{-\eta\kappa(t-\tau)} \vert z_r(\tau) \vert \label{eq: set point - intermediate estimates for lamda_k wk - 2}
\end{align}
for some constant $C_6 > 0$. Finally, we also have, for $-\gamma < -\tilde{\kappa} < - \kappa < 0$,
\begin{align}
\left\vert \int_0^t e^{\lambda_k(t-s)} \dot{r}_k(s) \,\mathrm{d}s \right\vert
& \leq \int_0^t e^{-\gamma(t-s)} \vert \dot{r}_k(s) \vert \,\mathrm{d}s \nonumber \\
& \leq \int_0^t e^{-(\gamma-\tilde{\kappa})(t-s)} \times e^{-\tilde{\kappa}(t-s)} \vert \dot{r}_k(s) \vert \,\mathrm{d}s \nonumber \\
& \leq \sqrt{ \int_0^t e^{-2(\gamma-\tilde{\kappa})(t-s)} \,\mathrm{d}s } \sqrt{ \int_0^t e^{-2\tilde{\kappa}(t-s)} \vert \dot{r}_k(s) \vert^2 \,\mathrm{d}s } \nonumber \\
& \leq \sqrt{ \dfrac{1}{2(\gamma-\tilde{\kappa})}} \sqrt{ \int_0^t e^{-2\tilde{\kappa}(t-s)} \vert \dot{r}_k(s) \vert^2 \,\mathrm{d}s } .
\end{align}
\end{subequations}
Taking the square on both sides of (\ref{eq: lambda_k w_k}), using Young's inequality, substituting estimates (\ref{eq: stability result - local exp stab}), (\ref{eq: prel estimate R}), and (\ref{eq: set point - intermediate estimates for lamda_k wk}), and using the fact that $\sum_{\vert k \vert \geq N+1} \vert \left< z , f_k \right> \vert^2 \leq \Vert z \Vert_\mathcal{H}^2/m_R$ for all $z \in \mathcal{H}$, we infer the existence of a constant $C_7 > 0$ such that
\begin{align*}
\sum\limits_{\vert k \vert \geq N+1} \vert \lambda_k w_k(t) \vert^2
& \leq C_7 e^{-2 \kappa t} \sum\limits_{\vert k \vert \geq N+1} \vert \lambda_k w_k(0) \vert^2 \\
& \phantom{\leq}\, + C_7 e^{-2\kappa t} \mathrm{CI}^2 + C_7 \sup\limits_{0 \leq \tau \leq t} e^{-2\eta\kappa(t-\tau)} \vert z_r(\tau) \vert^2 \\
& \phantom{\leq}\, + C_7 \int_0^t e^{-2\tilde{\kappa}(t-s)} \left\Vert \dfrac{\mathrm{d}R}{\mathrm{d}t}(s,\cdot) \right\Vert_\mathcal{H}^2 \,\mathrm{d}s ,
\end{align*}
where we recall that $\mathrm{Re}\lambda_k \leq -\gamma < -\kappa < 0$ for all $\vert k \vert \geq N+1$. Noting that, for $\vert k \vert \geq N+1 \geq n_0+1$, $\left< \mathcal{A}W(0) , f_k \right>_\mathcal{H} = \lambda_k w_k(0)$, we infer that $\sum_{\vert k \vert \geq N+1} \vert \lambda_k w_k(0) \vert^2 \leq \sum_{k \in \Z} \vert \left< \mathcal{A}W(0) , f_k \right>_\mathcal{H} \vert^2 \leq \Vert \mathcal{A} W(0) \Vert_\mathcal{H}^2 / m_R$. Then, based on (\ref{eq: stability result - local exp stab}) and (\ref{eq: estimate norm R_t}), and as $\tilde{\kappa} > \kappa > 0$, estimations similar to the ones reported in (\ref{eq: set point - intermediate estimates for lamda_k wk - 1}-\ref{eq: set point - intermediate estimates for lamda_k wk - 2}) show the existence of a constant $C_8 > 0$ such that 
\begin{align*}
\sum\limits_{\vert k \vert \geq N+1} \vert \lambda_k w_k(t) \vert^2
& \leq C_8 e^{-2\kappa t} \left\{ \mathrm{CI}^2 + \Vert \mathcal{A} W(0) \Vert_\mathcal{H}^2 \right\} + C_8 \sup\limits_{0 \leq \tau \leq t} e^{-2\eta\kappa(t-\tau)} \vert z_r(\tau) \vert^2 .
\end{align*}
Substituting this latter estimate and (\ref{eq: stability result - local exp stab}) into (\ref{eq: prel estimate tracking performance}), we obtain the existence of constants $\overline{C}_3, \overline{C}_4 > 0$ such that (\ref{eq: thm ref tracking - estimate}) holds. 
\qed

\section{Numerical illustration}\label{sec: numerical illustration}

For numerical illustration, we set $L=1$, $\alpha = 1.1$, and we consider the nonlinear function $f(y) = y^3$, which leads to the boundary control system:
\begin{subequations}\label{eq: numerical illustration - wave equation}
\begin{align}
&\dfrac{\partial^2 y}{\partial t^2} = \dfrac{\partial^2 y}{\partial x^2} + y^3 ,  \\
& y(t,0) = 0 , \qquad \dfrac{\partial y}{ \partial x}(t,L) = u(t) ,\\
& y(0,x) = y_0(x) , \qquad \dfrac{\partial y}{\partial t}(0,x) = y_1(x) , 
\end{align}
\end{subequations}
for $t > 0$ and $x \in (0,1)$. Since $F(y) = \int_0^y f(s) \,\mathrm{d}s = y^4/4 \rightarrow + \infty$ when $\vert y \vert \rightarrow + \infty$, it follows from Remark~\ref{rem: existence of steady states} the existence of a steady state $y_e \in \mathcal{C}^2([0,L])$ associated with any given value of the system output $z_e \in \R = y_e'(0)$. We set $z_e = 1.5$ and numerically compute the associated steady-state trajectory $y_e$, giving in particular the equilibrium control input $u_e = y_e'(L) \approx 0.781$. In the absence of non linearity, i.e. $f = 0$, it is known (see e.g.~\cite[Section~4]{russell1978controllability}) that the eigenvalues of the operator $\mathcal{A}$ defined by (\ref{eq: def operator A}) are given by
\begin{equation*}
\lambda_k = \dfrac{1}{2L} \log\left(\dfrac{\alpha-1}{\alpha+1}\right) + i \dfrac{k\pi}{L} .
\end{equation*}
These values are used as initial guesses to determine the eigenvalues $\lambda_k$ and the associated eigenvectors $e_k$ of the operator $\mathcal{A}$ in the presence of the non linearity $f(y)=y^3$ by using a shooting method. We obtain one unstable eigenvalue $\lambda_0 \approx 0.326$ while all other eigenvalues are stable with a real part less than $1$. The feedback gain is computed to place the poles of the truncated model at $-0.5$, $-1$, and $-1.5$.

\begin{figure}
\centering
	\includegraphics[width=3in]{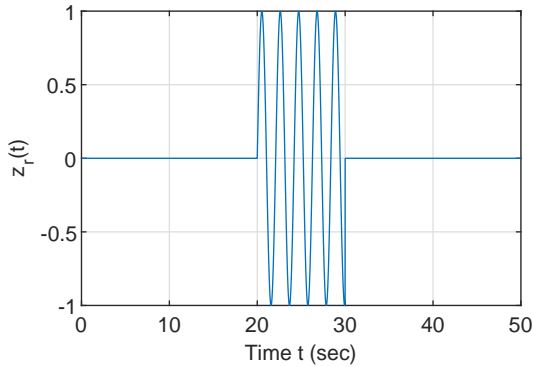}
	\caption{Signal $z_r(t)$}
	\label{fig: Ref_zr}	
\end{figure}

\begin{figure}
\centering
	\subfigure[State $y(t,x)$]{
	\includegraphics[width=2.8in]{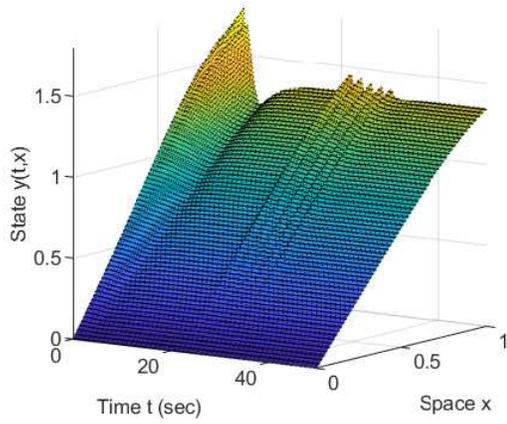}
	}
	\subfigure[State $\frac{\partial y}{\partial t}(t,x)$]{
	\includegraphics[width=2.8in]{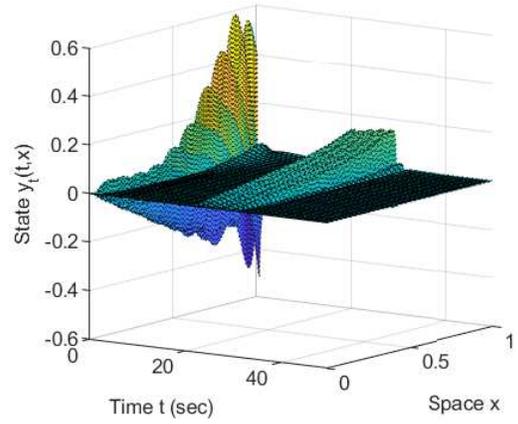}
	}
	\subfigure[Regulated output $z(t) = \frac{\partial y}{\partial x}(t,0)$]{
	\includegraphics[width=2.8in]{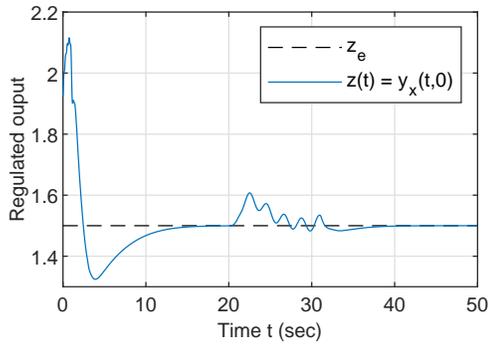}
	}
	\subfigure[Control input $u(t) = \frac{\partial y}{\partial x}(t,L)$]{
	\includegraphics[width=2.8in]{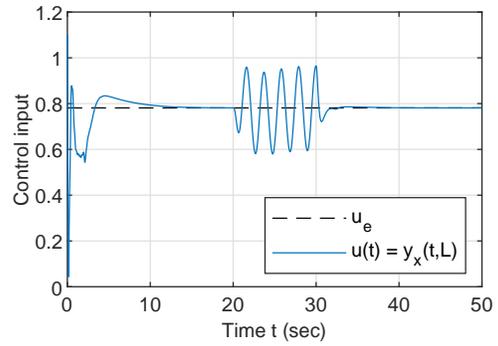}
	}
\caption{Time domain evolution of the state of the closed-loop system}
\label{fig: sim closed-loop system}
\end{figure}

For numerical simulations, we select the initial condition $W(0,x) = (\frac{2\alpha}{5}x,-\frac{2}{5L}x) \in D(\mathcal{A})$ and the signal $z_r$ as depicted in Fig.~\ref{fig: Ref_zr}. The adopted numerical scheme is the modal approximation of the infinite-dimensional system using its first 10 modes. The time domain evolution of the state of the closed-loop system, the regulated output, and the command input, are depicted in Fig.~\ref{fig: sim closed-loop system}. The simulation results are compliant with the theoretical predictions.

\section{Conclusion and open issues}\label{sec: conclusion}
In this paper we have investigated the proportional integral (PI) boundary regulation control of the left Neumann trace of a one-dimensional semilinear wave equation. Our control strategy combines a traditional velocity feedback and the design of an auxiliary control law performed on a finite-dimensional (spectrally) truncated model. 

A number of open issues and potential directions for further research emerge from this study, that we list and comment hereafter.

\subsection{Other controls and outputs}
\paragraph{Other controls}
In this paper, we have considered a Neumann boundary control, on the right of the interval. Of course, we could have taken this control on the left, or even on both sides. While the proposed PI procedure seems to be extendable to the cases of a Robin boundary control, of mixed Dirichlet-Neumann boundary controls, and of a distributed (i.e., internal) control input (note that, in this case, the control operator is bounded), the case of a Dirichlet boundary control seems much more challenging and remains open.

\paragraph{Other outputs}
We have selected the regulated ouput as the Neumann trace. However, one might be interested in regulating other types of outputs. This includes, for example, the Dirichlet trace at the left boundary or the value of $y$ either at a specific location or on a subset of the spatial domain. It is not clear whether or not the reported PI control design procedure could be extended to these different settings.

\subsection{The general multi-dimensional case}\label{sec_gen_multiD}
In this paper we have focused on the one-dimensional wave equation, which is already quite challenging. The multi-dimensional case is completely open, in particular because, except in very particular cases (like the two-dimensional disk for instance), we do not have any Riesz basis formed by generalized eigenvectors of the underlying unbounded operator. Indeed, the Riesz basis property is essentially a one-dimensional property.

Hence, in the multi-dimensional case, a first challenging difficulty that must be addressed is to replace the Riesz basis study by more abstract projection operators. On this issue, it seems that the recent work \cite{takahashi} provides very interesting line of research, although the authors of that reference deal with a parabolic equation with a selfadjoint operator and then one has a spectral decomposition and spectral projectors. For a wave equation, the underlying operator is not selfadjoint and is not even normal, which creates deep difficulties for the spectral study.

At present, we leave this issue completely open.

\subsection{Robustness with respect to perturbations}
Beyond the set-point regulation control of the system output, PI controllers are well known for their ability to provide, in general, a form of robustness with respect to external disturbances. The case of a continuously differentiable additive distributed disturbance $d(t,\cdot) \in L^2(0,L)$, yielding the dynamics
\begin{equation*}
\dfrac{\partial^2 y}{\partial t^2} = \dfrac{\partial^2 y}{\partial x^2} + f(y) + d ,
\end{equation*}
can be easily handled in our approach by merely embedding the contribution of the disturbance $d$ into the residual term $r$ given by (\ref{eq: def residual term}). See also~\cite{lhachemi2019pi} for the case of a linear reaction-diffusion equation. A similar robustness issue would be to evaluate the impact of an additive boundary disturbance $p(t)$ applying to the control input:
\begin{equation*}
\dfrac{\partial y}{ \partial x}(t,L) = u(t) + p(t) .
\end{equation*}
However, comparing to the case of distributed perturbations, the study of the robustness of infinite-dimensional systems with respect to boundary perturbations is generally much more challenging~\cite{mironchenko2019input}. 

\subsection{Robustness with respect to the domain}
Another interesting robustness issue regarding the studied wave equation (\ref{eq: wave equation}) concerns the robustness of the PI procedure with respect to uncertainties on the length $L > 0$ of the domain. Indeed, in many engineering devices the value of $L$ is known only approximately; it may even happen that, along the process, the domain varies a bit. The robustness of the control strategy with respect to $L$ becomes then a major issue.

In our strategy, the value of the parameter $L$ directly impacts the vectors of the Riesz-basis $(e_k)_{k \in\Z}$. However, since the control strategy relies on a finite dimensional truncated model, only the perturbations on a finite number of vectors of the Riesz-basis have an impact on the control strategy. Hence, similarly to the study reported in~\cite{coron2006global}, the Riesz basis property remains uniform with respect to small perturbations of $L > 0$. This allows to easily derive a robustness result versus small variations of the length of the spatial domain, and thus to give a positive answer to the question of robustness with respect to $L$.

This result shows the power and advantage of the strategy developed here, based on spectral finite-dimensional truncations of the problem, which conveys robustness properties.

In the multi-dimensional case, the question of robustness with respect to the domain would certainly be much more difficult to address. It would be required to introduce an appropriate shape topology and then perform a kind of sensitivity analysis on the control design with respect to that topology.

\subsection{Robustness with respect to input and/or state delays}
In this section, we discuss the impact of possible input or state delays on the reported control strategy. As we explain below, we cannot expect any robustness of our strategy with respect to input delays, while the situation for state delays is expected to be more favourable.

\paragraph{Input delay}
It was shown in~\cite{lhachemi2019pi} for a 1-D reaction-diffusion equation that the PI controller design procedure can be augmented with a predictor component~\cite{artstein1982linear,bresch2018new,lhachemi2018feedback} to handle constant control input delays. Such a strategy fails in the context of the wave equation studied in this paper. The reason is that, even for arbitrarily small $h>0$, the preliminary feedback
\begin{equation*}
u_\delta(t) = - \alpha \dfrac{\partial y_\delta}{ \partial t}(t-h,L)
\end{equation*}
might fail to stabilize the linear wave equation, yielding an infinite number of unstable modes. To illustrate this remark, let us consider the following linear wave equation:
\begin{subequations}\label{eq: wave equation - input delay}
\begin{align}
& \dfrac{\partial^2 y}{\partial t^2} = \dfrac{\partial^2 y}{\partial x^2} - \beta y , \\
& y(t,0) = 0 , \qquad \dfrac{\partial y}{ \partial x}(t,L) = - \alpha \dfrac{\partial y}{\partial t}(t-h,L) ,
\end{align}
\end{subequations}
where $h > 0$ is a constant delay, $\alpha > 0$, and $\beta \in \R$. It was shown in~\cite{datko1997two} in the case $\alpha = 1$ and $\beta = 0$ that (\ref{eq: wave equation - input delay}) admits an infinite number of unstable modes for some arbitrarily small values of the delay $h > 0$. We elaborate here this remark by extending it to the case $\alpha > 0$ and $\beta\in\R$. We proceed similarly to~\cite{datko1997two} by looking for solutions of (\ref{eq: wave equation - input delay}) under the form $w(t,x) = e^{\lambda t} \sinh\left( \sqrt{\lambda^2+\beta} x \right)$ for some $\lambda\in\C$. It is easy to see that such a solution exists if and only if $\lambda\in\C$ satisfies
\begin{equation*}
\sqrt{\lambda^2+\beta} \cosh\left( \sqrt{\lambda^2 + \beta}L \right) = -\alpha\lambda e^{-\lambda h} \sinh\left( \sqrt{\lambda^2 + \beta}L \right) .
\end{equation*}
We start by studying the special case $\beta = 0$. Then the above identity becomes equivalent to:
\begin{equation}\label{eq: wave delay - case beta = 0}
e^{\lambda h} = - \alpha \tanh(\lambda L) .
\end{equation}
Let $k \in \Z$ be arbitrarily fixed and select the input delay
\begin{equation}\label{eq: wave delay - def h}
h = \dfrac{L}{k+1/2} .
\end{equation}
Let $\gamma > 0$ be the only positive number such that
\begin{equation}\label{eq: wave delay - def gamma}
e^{\gamma h} = \alpha \coth(\gamma L) .
\end{equation}
It is easily seen that $\lambda^0_n = \gamma + \dfrac{i}{L} \left( k + \dfrac{1}{2} \right) (4n+1)\pi$, $n\in\Z$, 
are distinct solutions of (\ref{eq: wave delay - case beta = 0}) with $\operatorname{Re}\lambda_n = \gamma > 0$. We now turn out attention onto the case $\beta \neq 0$. Let $h,\gamma > 0$ be still given by (\ref{eq: wave delay - def h}-\ref{eq: wave delay - def gamma}). We introduce the open set
$A = \left\{ \lambda\in \C \,:\, 0 < \operatorname{Re}z < 2 \gamma \;, \vert\lambda\vert > \sqrt{\beta} \right\}$.
We define for $\lambda \in A$ the holomorphic functions
\begin{align*}
f(\lambda) & = \lambda \cosh\left( \lambda L \right) + \alpha\lambda e^{-\lambda h} \sinh\left( \lambda L \right) , \\
g(\lambda) & = \sqrt{\lambda^2+\beta} \cosh\left( \sqrt{\lambda^2 + \beta}L \right) + \alpha\lambda e^{-\lambda h} \sinh\left( \sqrt{\lambda^2 + \beta}L \right) ,
\end{align*}
with 
$\sqrt{\lambda^2+\beta} = \lambda \exp\left(\dfrac{1}{2}\operatorname{Log}\left( 1 + \dfrac{\beta}{\lambda^2} \right)\right)$, 
where $\operatorname{Log}$ denotes the principal determination of logarithm. We have $f(\lambda_n^0) = 0$ for all $n \in \Z$ and it can be observed that 
\begin{equation}\label{eq: wave delay - Rouche thm - step 1}
g(\lambda) = f(\lambda) + O(1) , \qquad \lambda\in A , \quad \vert \lambda \vert \rightarrow + \infty .
\end{equation}
For a given constant $R > 0$ to be defined later, we consider the simple loop $\lambda_n : [0,2\pi] \rightarrow \C$ defined by
$\lambda_n(\theta) = \lambda_n^0 + \dfrac{R e^{i\theta}}{n}$, $\theta \in [0,2\pi]$.
We consider an integer $N_0 \geq 1$ such that $\lambda_n(\theta) \in A$ for all $\vert n \vert \geq N_0$ and all $\theta \in [0,2\pi]$. Standard computations show that
\begin{equation*}
f(\lambda_n(\theta)) = (-1)^k i \lambda_n^0 \dfrac{Re^{i\theta}}{n} \zeta \cosh(\gamma L) + O\left(\dfrac{1}{n}\right)
\end{equation*}
when $\vert n \vert \rightarrow + \infty$, uniformly with respect to $\theta\in[0,2\pi]$, where
$\zeta = L + h\tanh(\gamma L) - L \tanh^2(\gamma L)$.
We note that $\zeta \neq 0$. Indeed, the condition $\zeta=0$ would imply that $\tanh(\gamma L) > 0$ is the positive root of $L + h X - L X^2$, yielding the contradiction $\tanh(\gamma L) = (h+\sqrt{h^2+4L^2})/(2L) > 1 + h/(2L) > 1$. We deduce that
\begin{equation}\label{eq: wave delay - Rouche thm - step 2}
\vert f(\lambda_n(\theta)) \vert \sim C R
\end{equation}
when $\vert n \vert \rightarrow + \infty$, uniformly with respect to $\theta\in[0,2\pi]$, where
the constant $C = \dfrac{2(2k+1)\pi \vert \zeta \vert \cosh\left(\gamma L\right)}{L} \neq 0$
is independent of $R > 0$. Hence, in view of (\ref{eq: wave delay - Rouche thm - step 1}-\ref{eq: wave delay - Rouche thm - step 2}) and by selecting $R > 0$ large enough, we can apply Rouch{\'e}'s theorem for $\vert n \vert \geq N_1$ with $N_1 > 0$ large enough. This shows for any $\vert n \vert \geq N_1$ the existence of
$\lambda_n^\beta = \lambda_n^0 + O\left(\dfrac{1}{n}\right)$
such that $g(\lambda_n^\beta) = 0$. Since $h > 0$ defined by (\ref{eq: wave delay - def h}) can be made arbitrarily small by selecting $k$ arbitrarily large, we have shown the existence of arbitrarily small delays $h > 0$ such that (\ref{eq: wave equation - input delay}) admits an infinite number of unstable modes. Such an observation implies that the strategy reported in this paper cannot be successfully applied to the case of a delayed control input. Hence, the PI regulation control of (\ref{eq: wave equation}) in the presence of a delay in the control input remains open.

\paragraph{State delay}
While, as discussed above, the case of a delay in the boundary control induces many difficulties (in particular, an infinite number of unstable modes), the case of a delay in the non-linearity $f$ might be more favourable. More specifically, a possible research direction is concerned with the potential extension of the control strategy reported in this paper to the case of the wave equation
\begin{equation*}
\dfrac{\partial^2 y}{\partial t^2}(t,x) = \dfrac{\partial^2 y}{\partial x^2}(t,x) + f(y(t-h,x)) ,
\end{equation*}
where $h > 0$ is a state-delay. Possible approaches include the use of either Lyapunov-Krasovskii functionals~\cite{kolmanovskii2012applied} or small gain arguments~\cite{lhachemi2020boundary}.

\subsection{Alternatives for control design}
In this paper, we have designed the control strategy, in particular, by using the pole shifting theorem (see in particular equation~\ref{eq: control input vd} in Subsection \ref{sec_poleshift}). Such a pole shifting is very natural considering the spectral approach that we have implemented: placing adequately the poles (of course, the unstable ones; but we can also shift to the left some of the stable ones in order to improve the stabilization properties) we are able to ensure some robustness properties, with respect to disturbances or with respect to the domain, as we have discussed previously.

One can wonder whether this is possible to design the control by other methods. For instance by the celebrated Riccati procedure (see, e.g., \cite{Khalil,Sontag,trelatbook,Zabczyk}) applied to the truncated finite-dimensional system. Considering the classical Linear Quadratic Riccati theory, let us denote by $u_n$ the optimal Riccati control (for some given weights) associated with the truncated system in dimension $n$. We do not make precise all notations but the framework is clear. The main question is: does $u_n$ converge to $u_\infty$ as $n\rightarrow+\infty$, where $u_\infty$ is the Riccati control of the complete system in infinite dimension (i.e., the PDE)?

We expect that such a convergence property is true at least in the parabolic case, i.e., for instance, for heat-like equations, or in the hyperbolic case with internal controls, i.e., for instance, for wave-like equations with distributed control. The hyperbolic case with boundary controls is probably much more difficult.

Adjacently, this discussion raises the problem of discussing the numerical efficiency of various control design approaches. It would be interesting to compare our approach, developed in the present paper, with other possible approaches and to compare their efficiency. One different approach that may come to our mind is backtepping design, which is well known to promote robustness properties (see, e.g., \cite{Krstic}), at least, for parabolic equations. It is not clear whether backstepping could be performed for the 1-D wave equation investigated in the present article.

\subsection{Controlling the output regulation by quasi-static deformations}

In the present work, we have dealt with the stabilization and regulation control in the vicinity of a given steady-state. 
In this section, we address the following question: considering that the system output has been regulated to a given steady-state, is it now possible to steer the system output to \emph{another} steady-state?

More precisely, based on Definition \ref{def: steady state}, we denote by $\mathcal{S} \subset\mathcal{C}^2([0,L])$ the set of steady-states of (\ref{eq: wave equation}) endowed with the $\mathcal{C}^2([0,L])$ topology. Let $y_{e,1},y_{e,2} \in \mathcal{S}$ belonging to the same connected component of $\mathcal{S}$. 
Introducing for $i \in \{1,2\}$ the system output associated with the steady-state $y_{e,i}$ defined by
\begin{equation*}
z_{e,i} = \dfrac{\mathrm{d} y_{e,i}}{\mathrm{d} x}(0),
\end{equation*}
can we design a PI controller 
able to steer the system output from its initial value $z_{e,1}$ (or, close to it) into any neighbourhood of $z_{e,2}$ in finite time? 

\medskip

A way to address this issue is to design a PI controller by \emph{quasi-static deformation}, as in \cite{coron2004global,coron2006global,schmidt2006}, along a path of steady outputs connecting $z_{e,1}$ to $z_{e,2}$. More precisely, since $y_{e,1},y_{e,2} \in \mathcal{S}$ are assumed to belong to the same connected component of $\mathcal{S}$, one can define $y_e(\tau,\cdot)$, with $\tau \in [0,1]$, a $\mathcal{C}^2$ path in $\mathcal{S}$ connecting $y_{e,1}$ to $y_{e,2}$, as well as the associated path of boundary inputs $u_e(\tau) = \frac{\partial y_e}{\partial x}(\tau,L)$.
Now, a possible attempt would be to extend the approach of~\cite{coron2006global} by defining the path of system outputs $z_e(\tau) = \frac{\partial y_e}{\partial x}(\tau,0)$ associated with $y_e(\tau,\cdot)$. In this setting, the objective would be to study the deviations of the system output $z(t) = \frac{\partial y}{ \partial x}(t,0)$ with respect to the quasi-static path $y_e(\epsilon t , \cdot)$ by introducing $z_\delta(t,x) =  z(t) - z_e(\epsilon t)$.
In this configuration, $r(t) = z_e(\epsilon t)$ plays the role of a slowly time-varying reference input.

 It was shown in~\cite{coron2006global} that such a path can be used in order to steer the system (\ref{eq: wave equation}) from the steady-state $y_{e,1}$ to the steady-state $y_{e,2}$ by means of a boundary control input $u$ taking advantage of quasi-static deformations. Specifically, for $\epsilon > 0$ small enough, the authors studied the deviations of the system trajectory with respect to the quasi-static path $y_e(\epsilon t , \cdot)$ by introducing:
\begin{align*}
y_\delta(t,x) & =  y(t,x) - y_e(\epsilon t , x) , \\
u_\delta(t) & = u(t) - u_e(\epsilon t) 
\end{align*}
for $t \in [0,1/\epsilon]$ and $x \in [0,L]$. The preliminary feedback still takes the form (\ref{eq: preliminary control input}). Due to the quasi-static deformations-based approach and using a Taylor expansion as in (\ref{eq: wave equation - variations around equilibrium - PDE}), the design of the auxiliary control input $v(t)$ requires the introducing of the following family of wave operators parametrized by $\tau \in [0,1]$:
\begin{equation*}
\mathcal{A}(\tau) = 
\begin{pmatrix}
0 & \mathrm{Id} \\ \mathcal{A}_0(\tau) & 0
\end{pmatrix}
\end{equation*}
with $\mathcal{A}_0(\tau) = \Delta + f'(y_e(\tau,\cdot)) \,\mathrm{Id}$ defined on the domain
\begin{align*}
D(\mathcal{A}(\tau)) 
= \{ (w^1,w^2) \in \mathcal{H} \,:\, & w^1 \in H^2(0,L) ,\, w^2 \in H^1(0,L) ,\, \\
& w^2(0) = 0 ,\, (w^1)'(L)+\alpha w^2(L) = 0 \} .
\end{align*} 
Following~\cite[Lem.~2]{coron2006global}, this family of operators admits a family $(e_{k}(\tau,\cdot))_{k \in \Z}$ of Riesz-bases formed by generalized eigenvectors of $\mathcal{A}(\tau)$, associated to the eigenvalues $(\lambda_k(\tau))_{k \in \Z}$ and with dual Riesz basis $(f_k(\tau,\cdot))_{k \in \Z}$, with properties similar to the ones of (\ref{lem: properties A}) but with an integer $n_0 \geq 0$ that is uniform with respect to $\tau \in [0,1]$. Without loss of generality, this latter integer can be selected such that $\vert k \vert \geq n_0 + 1$ implies $\operatorname{Re}\lambda_k < -1$. Moreover, the aforementioned family of Riesz-bases is uniform with respect to $\tau \in [0,1]$ in the sense that the constants $m_R,M_R$ of Definition \ref{def: Riesz basis} can be selected independently of $\tau \in [0,1]$. These key properties allowed the authors of~\cite{coron2006global} to design a control law of the form $v(t) = K(\epsilon t) X(t)$. The matrix $K(\tau)$, parametrized by $\tau \in [0,1]$, is obtained based on an augmented finite-dimensional LTI system, also parametrized by $\tau \in [0,1]$, which in particular captures the first modes of $\mathcal{A}(\tau)$ characterized by the integers $- n_0 \leq k \leq n_0$. The vector $X(t)$ captures, in addition to a number of integral components, the projection of the system trajectory $W(t,\cdot)$ onto the vector space spanned by $(e_{k}(t\epsilon,\cdot))_{\vert k \vert \leq n_0}$. The stability property of the resulting closed-loop system was assessed through the study of a suitable Lyapunov functional, yielding the following result~\cite[Thm.~1]{coron2006global}. For the wave equation (\ref{eq: wave equation}) with initial condition set as the steady-state $y_{e,1}$ and for the control input selected as above: for every $\delta > 0$, there exists $\epsilon_1 > 0$ so that, for every $\epsilon \in (0,\epsilon_1]$, we have
\begin{equation*}
\left\Vert \dfrac{\partial y}{\partial x}(1/\epsilon,\cdot) - \dfrac{\mathrm{d} y_{e,2}}{\mathrm{d}x} \right\Vert_{L^2(0,L)}
+ \left\Vert \dfrac{\partial y}{\partial t}(1/\epsilon,\cdot) \right\Vert_{L^2(0,L)}
\leq \delta .
\end{equation*} 
In conclusion, it is of interest to evaluate the possible extension of our PI regulation procedure to the case of quasi-steady deformations as described above.

\subsection{System of one-dimensional partial differential equations}
In Section \ref{sec_gen_multiD}, we have mentioned as a completely open issue the general multi-dimensional case, which is very challenging.
As an intermediate case between 1-D and multi-D, we may consider the case of systems of one-dimensional partial differential equations.
A line of research that is of great interest is to consider coupled scalar one-dimensional PDEs, for instance of the form
\begin{align*}
\dfrac{\partial^2 y}{\partial t^2} &= a \dfrac{\partial^2 y}{\partial x^2} + c y + dz \\
\dfrac{\partial^2 z}{\partial t^2} &= b \dfrac{\partial^2 y}{\partial x^2} + e y + fz 
\end{align*}
which are 1-D coupled wave-like equations, with various possible controls and with various outputs (e.g., Neumann boundary control and Neumann regulated output like in this article).

Of great interest too, would be to replace the above wave equation in $z$, with a heat-like equation in $z$, that is, consider a 1-D wave-like equation that is coupled with a 1-D heat-like equation in $z$. In that case, we expect new phenomena, emerging from the interesting coupling between a parabolic and a hyperbolic equation. Actually, even the case of coupled 1-D heat-like equations does not seem to have been considered in the literature concerning PI issues. 

It would be very interesting to address these problems. They indeed have attracted much attention for controllability issues and many powerful techniques have been introduced to treat such coupled systems (see, e.g., \cite{Alabau2002,Alabau2015,AlabauLeautaud,AK1,AK2,AK3,LiardLissy}).
This open issue is a future line of research.

\appendix

\section{Annex - Proof of Lemma~\ref{lem: properties A}}\label{annex: proof lemma}

From the definition of the operator $\mathcal{A}$ given by (\ref{eq: def operator A}), $\lambda \in\C$ is an eigenvalue of $\mathcal{A}$ associated with the nonzero eigenvector $w = (w^1,w^2) \in D(\mathcal{A})$ if and only if $w^2 = \lambda w^1$ and
\begin{align*}
& (w^1)'' + f'(y_e) w^1 = \lambda^2 w^1 , \\
& w^1(0)=0 , \quad (w^1)'(L)+\alpha\lambda w^1(L) = 0
\end{align*}
for $x \in (0,L)$. Then, for $\vert \lambda \vert \rightarrow + \infty$, we obtain for $x \in [0,L]$ that
\begin{equation*}
w^1(x) = \sinh\left(\sqrt{\lambda^2+O(1)} x\right) , \quad
(w^1)'(x) = \sqrt{\lambda^2+O(1)} \cosh\left(\sqrt{\lambda^2+O(1)} x\right) ,
\end{equation*}
uniformly with respect to $x \in [0,L]$. Using now the right boundary conditions, we obtain the existence of an integer $k \in\Z$ such that, as $\vert k \vert \rightarrow +\infty$, 
\begin{equation*}
\lambda_k = \dfrac{1}{2L} \log\left(\dfrac{\alpha-1}{\alpha+1}\right) + i \dfrac{k\pi}{L} + O\left( \dfrac{1}{\vert k \vert} \right) .
\end{equation*}
Then, an associated unit eigenvector is given by
\begin{equation*}
e_k = \dfrac{1}{A_k} \left( \sinh\left(\sqrt{\lambda_k^2+O(1)} x\right) , \lambda_k \sinh\left(\sqrt{\lambda_k^2+O(1)} x\right) \right)
\end{equation*}
where, recalling that $\alpha > 1$ and introducing $\beta = - \frac{1}{2L} \log\left(\frac{\alpha-1}{\alpha+1}\right) > 0$, 
\begin{equation*}
A_k = \vert \lambda_k \vert \sqrt{ \dfrac{\sinh(2\beta L)}{2\beta} + O\left(\dfrac{1}{\vert k \vert}\right) } .
\end{equation*}
In particular, $(e_k^1)'(0) = O(1)$ as $\vert k \vert \rightarrow + \infty$, showing item 7.

We show that the eigenvalues of $\mathcal{A}$ are geometrically simple (item 2). To do so, assume that $w_i = (w_i^1,w_i^2) \in D(\mathcal{A})$, with $i\in\{1,2\}$, are two eigenvectors of $\mathcal{A}$ associated with the same eigenvalue $\lambda \in\C$. We note that $w_i(L) \neq 0$ because otherwise $(w_i^1)'(L) = - \alpha \lambda w_i^1(l) = 0$ hence, by cauchy uniqueness, $w_i^1 = 0$ and thus $w_i^2 = \lambda w_i^1 = 0$, giving the contradiction $w = 0$. Then the function $g$ defined by $g = w_2^1(L) w_1^1 - w_1^1(L) w_2^1 \neq 0$ satisfies
\begin{align*}
& g'' + f'(y_e) g = \lambda^2 g , \\
& g(L)=g'(L)=0
\end{align*}
implying $g=0$. Recalling that $w_i^2 = \lambda w_i^1$, this shows that $w_1$ and $w_2$ are not linearly independent.

Recalling that $\mathcal{A}$ has compact resolvent, we denote by $(e_k)_{k\in\Z}$ a complete set of unit generalized eigenfunctions of $\mathcal{A}$ associated with the eigenvalues $(\lambda_k)_{k\in\Z}$~\cite{gohberg1978introduction}. We are going to apply Bari's theorem~\cite{gohberg1978introduction} to show that $(e_k)_{k\in\Z}$ is a Riesz basis. To do so, we need a Riesz basis of reference. Based on the definition the operator $\mathcal{A}$ given by (\ref{eq: def operator A}), we consider the below operator, obtained by removing the contribution of the terme $f'(y_e) \,\mathrm{Id}$,
\begin{equation*}
\mathcal{A}_{tr} = \begin{pmatrix}
0 & \mathrm{Id} \\ \Delta & 0
\end{pmatrix}
\end{equation*}
defined on the same domain as $\mathcal{A}$. We know from~\cite[Section~4]{russell1978controllability} that $\mathcal{A}_{tr}$ admits a Riesz basis of eigenvectors $(\phi_k)_{k\in\Z}$ associated with the eigenvalues $(\mu_k)_{k\in\Z}$ given for any $k\in\Z$ by 
\begin{equation*}
\mu_k = \dfrac{1}{2L} \log\left(\dfrac{\alpha-1}{\alpha+1}\right) + i \dfrac{k\pi}{L} ,
\end{equation*}
and 
\begin{equation*}
\phi_k = \dfrac{1}{B_k} \left( \sinh(\mu_k x) , \mu_k \sinh(\mu_k x) \right)
\end{equation*}
where, recalling that $\beta = - \frac{1}{2L} \log\left(\frac{\alpha-1}{\alpha+1}\right) > 0$, 
\begin{equation*}
B_k = \dfrac{1}{L\sqrt{2\beta}} \sqrt{(\beta^2 L^2 + k^2 \pi^2)\sinh(2\beta L)} .
\end{equation*}
We deduce that
\begin{equation*}
e_k = \phi_k + O\left(\dfrac{1}{\vert k \vert}\right) ,
\end{equation*}
in $\mathcal{H}$-norm as $\vert k \vert \rightarrow +\infty$. Hence $(e_k)_{k\in\Z}$ is quadratically close to the Riesz basis $(\phi_k)_{k\in\Z}$. Then, with the results of~\cite[Lemma~6.2 and Theorem~6.3]{guo2001riesz} relying on Bari's theorem, we obtain that $(e_k)_{k\in\Z}$ is a Riesz basis. Introducing $(f_k)_{k\in\Z}$ as the Dual Riesz basis of $(e_k)_{k\in\Z}$, items 1, 3, 4, and 6 hold true. Finally, an homotopy argument using the operator $\mathcal{A}_{tr}$ shows that the algebraic multiplicity of the real eigenvalues of $\mathcal{A}$ is odd, yielding item 5.

\bibliographystyle{IEEEtranS} 
\bibliography{mybibfile}

\end{document}